# Hamilton–Jacobi–Bellman–Isaacs equation for rational inattention in the long-run management of river environments under uncertainty

**Authors** Hidekazu Yoshioka[1, 2, *], Motoh Tsujimura[3]

[1] Graduate School of Natural Science and Technology, Shimane University, Nishikawatsu-cho 1060, Matsue, 690-8504, Japan

[2] Fisheries Ecosystem Project Center, Shimane University, Nishikawatsu-cho 1060, Matsue, 690-8504, Japan

[3] Graduate School of Commerce, Doshisha University, Karasuma-Higashi-iru, Imadegawa-dori, Kamigyo-ku, Kyoto, 602-8580, Japan

* Corresponding author
  E-mail: yoshih@life.shimane-u.ac.jp

**Abstract**

A new stochastic control model for the long-run environmental management of rivers is mathematically and numerically analyzed, focusing on a modern sediment replenishment problem with unique nonsmooth and nonlinear properties. Rational inattention as a novel adaptive strategy to collect information and intervene against the target system is modeled using Erlangization. The system dynamics containing the river discharge following a continuous-state branching with an immigration-type process and the controlled sediment storage dynamics lead to a nonsmooth and nonlocal infinitesimal generator. Modeling uncertainty, which is ubiquitous in certain applications, is considered in a robust control framework in which deviations between the benchmark and distorted models are penalized through relative entropy. The partial integro-differential Hamilton–Jacobi–Bellman–Isaacs (HJBI) equation as an optimality equation is derived, and its uniqueness, existence, and optimality are discussed. A monotone finite difference scheme guaranteeing the boundedness and uniqueness of numerical solutions is proposed to discretize the HJBI equation and is verified based on manufactured solutions. Model applications are also conducted with the parameter values identified from the available data and physical formulae. The computational results suggest that environmental management should be rationally inattentive in a state-dependent and adaptive manner.

## 1. Introduction

### 1.1 Modeling environmental management

Finding sustainable and sound management strategies for the environment is an important mission for society [1]. River environments play a crucial role in human life by providing fishery resources [2] and water resources for hydropower [3], irrigation [4], and drinking water [5], as well as in mitigating flood disasters [6]. Furthermore, river environments serve not only as hydrological pathways of global cycles of water, dissolved matter, and sediment [7], but also as essential aquatic habitats [8]. By contrast, human interventions against rivers can trigger severe environmental problems, including pollution [9], hydrological regime shifts and associated changes in benthic fauna [10, 11], and sediment starvation [12, 13]. Therefore, balancing between human intervention and environmental restoration is a critical issue.

Environmental dynamics in rivers are driven by stream flows, which are inherently stochastic owing to stochastic and unpredictable global hydrological dynamics [14, 15]. It is therefore natural to consider the problem of river management as a stochastic control problem [16]. Target system dynamics, including streamflow dynamics [17, 18], associated chemicals [19], and biological [20] and physical dynamics [21], can be formulated using stochastic differential equations (SDEs). With the help of an optimization principle such as dynamic programming [16], an optimal control can be obtained analytically or numerically.

### 1.2 Jump-driven streamflow dynamics

Streamflow dynamics are a key factor in the stochastic control of river environments. From an operational viewpoint, streamflow dynamics governing the temporal evolution of hydraulic quantities, such as discharges, should be simple and efficient. Such requirements can be achieved using SDEs driven by well-defined martingale noises such as Brownian noises [22], compound Poisson processes [23], Markov chains [17], and more general Lévy processes [18, 20].

As a common drawback of existing models, the noises have independent increments, suggesting that important but natural clustered events, such as rainfall-induced floods [24, 25], are not considered. Such events have state-dependent jumps with positive feedback and follow self-exciting jump processes [26, 27]. A physically more reasonable model with a satisfactory statistical performance is desired for better descriptions of the streamflow dynamics. Although generalized continuous-state branching processes with immigration (gCBI processes) [28] can be an efficient as well as well-defined candidate self-exciting jump processes, to the best of authors' knowledge, their application to environmental modeling has yet to be addressed.

### 1.3 Sediment replenishment

We focus on the sediment starvation issue in rivers as a crucial environmental problem encountered by many rivers with dams. Sediment replenishment is a countermeasure used to restore aquatic environments, such as dam-downstream rivers and coastal areas subject to severe erosion owing to the deficit of natural sediment supply. This is a human intervention in which sediment is transported from a source to a target place, posing an optimization problem to decide when and how much of the sediment should be replenished [29]. Fine materials, including sand, have been supplied to coastal areas [30] and relatively large materials such as gravel have been supplied to rivers [31]. Sand has also been replenished in dam-downstream rivers [29, 32]. For river restoration, sediment is replenished from once or twice [33] to up to 13 times [34] a year. In addition, different refreshment frequencies were used in different years [35].

The optimization of sediment replenishment has been considered from both static and dynamic perspectives. For example, McNamara et al. [36] considered the static problem of beach nourishment under Poisson storm arrivals. In addition, Cutler et al. [37] formulated a dynamic programming problem for beach nourishment under sea-level rise scenarios. A game-theoretic analytical approach was also employed by Janssen et al. [38]. Dynamic programming approaches of sediment replenishment in rivers where sediment storage dynamics follow SDEs have also been analyzed [29, 32, 39], although with extremely simplified streamflow dynamics using compound Poisson jumps.

### 1.4 Intervention models

Dynamic optimization of the amount and timing of interventions can be considered as an impulse control [16]. Conventional models assume that the decision-maker continuously observes the target dynamics the and can intervene at any time. This assumption applies if the target dynamics are frequently observable; however, this is not always true for applications in the environment and ecology because the observation is costly and difficult to carry out in a frequent manner. Weekly to monthly biological and environmental surveys of aquatic environments are possible [40, 41]. Environmental monitoring is also time-consuming [42]. Do and Thi [43] suggested that the optimization of the sampling frequency in natural environments under anthropogenic pressure should be adaptive. Sampling and intervention strategies should be optimized such that the observation cost becomes minimized while the dynamics are effectively controlled. This is a rational inattention problem [44] that balances the costs and benefits of controlling the dynamics by tuning the rate of information collection.

Several optimization models of stochastic dynamics based on rational inattention have been developed. The Poisson observation models used to observe the target dynamics at the jump times of a Poisson process are a major stream of mathematical models [45]. They are tractable and utilized in the finance and insurance fields [46, 47]. While the Poisson assumption provides simple Markovian models, it assumes that the sampling intervals follow an exponential distribution, which is possibly restrictive in terms of application. Semi-Markov models, as more flexible alternatives, allow for more flexible sampling intervals [48, 49]. Models with state-dependent deterministic observation intervals have also been discussed for queuing systems [50], contract design [51], disease treatment [52], and environmental management [53, 54].

As an efficient and flexible randomizing method, Avanzi et al. [55] considered semi-Markov dynamic programming of dividend decisions through Erlangization, bridging the random and deterministic sampling strategies. Erlangization is based on the probabilistic nature of the Erlang distributions, where the simplest case corresponds to the exponential distribution, whereas its appropriate scaling limit is a delta

(deterministic) distribution [56], serving as a randomization technique bridging the deterministic and stochastic sampling intervals. Erlangization has been applied to the field of insurance management [57, 58] and fluid queues [59, 60], but its applications to other problems are rare; in particular, we did not find applications in the field of environmental management as mentioned above.

**1.5 Model uncertainty**

A common issue in optimal control is the uncertainty caused by the deviation of a formulated model from a real model. Anderson et al. [61] established a dynamic programming approach to robustly control the stochastic dynamics and tolerate model uncertainty. In their framework, the uncertainty between the benchmark and distorted models is evaluated using the relative entropy measured through a Radon–Nikodým derivative. Therefore, it assumes an equivalence between the benchmark and distorted models. Despite this restriction, the robust control approach successfully handles complex problems, such as mean-variance insurance management [62], an economic risk assessment of climate changes [63], resources [64], and ecological management [65]. Its applications to self-exciting processes addressed in this paper, as explained below, are scarce.

**1.6 Objectives and contributions**

Our objectives are the formulation, analysis, and application of a stochastic control model for sustainable and sound development of river environments, where the decision-maker conducts discrete observations adaptively and the streamflow dynamics follow an uncertain gCBI process. This is a unique engineering problem arising in the modern world, which widely ranges from modeling and computations to the application of a partial integro-differential equation (PIDE) with nonsmooth coefficients and an exponential nonlinearity. Although in this study we focus on a specific problem, because (g)CBI processes have many engineering applications, the theoretical results obtained are widely applicable to other jump-driven systems with uncertainty.

We focus on the sediment replenishment problem as a model case and address four specific issues:

**(i) Modeling adaptive discrete observation/intervention:** Our problem is related to rational inattention in the sense that the opportunities of observing the target dynamics are optimized to balance the observation and other costs and penalties. Erlangization combined with impulse control realizes an adaptive management strategy to flexibly control the intensity of the observation. This introduces an innovative mathematical framework for partial observational control.

**(ii) Modeling through an uncertain gCBI process:** We incorporate the framework of Anderson et al. [61] into a gCBI process of streamflow dynamics. Sediment storage is considered a controlled process following an SDE with nonsmooth coefficients representing physical constraints. The jump kernel of the gCBI process is uncertain, and its distortion is expressed using a Radon–Nikodým derivative [66, 67]. We evaluate the relative entropy between the benchmark and distorted models, where the deviation between them is softly penalized in an objective functional. The problem then becomes finding the most cost-efficient observation/intervention strategy of the dynamics under uncertainty, i.e., a zero-sum stochastic differential game between the decision-maker as a minimizing player and nature as a maximizing opponent player. To date, CBI processes with uncertainties have not been considered for robust control.

(iii) **Solving the optimality equation:** Our control problem is reduced to solving the Hamilton–Jacobi–Bellman–Isaacs (HJBI) equation as a PIDE with nonsmooth coefficients. Although the problem does not admit closed-form solutions, the mathematical rigor of the viscosity solutions [68] applies, and we discuss whether the equation is optimal under certain regularity assumptions. We also provide the first mathematical results on the HJBI equation based on the gCBI process under uncertainty and discrete observation/intervention.

(iv) **Computation and Application:** We use a monotone finite difference scheme to compute the HJBI equation. The scheme guarantees the uniqueness and boundedness of numerical solutions, and its performance is verified based on manufactured solutions. Finally, our model is applied to a realistic case with coefficients and parameter values that are successfully identified using real data and semi-empirical formulae. We then compute the optimal observation/intervention strategies based on rational inattention and Erlangization.

Appendices contain proofs of the propositions (**Appendix A**) and several technical results on the gCBI processes (**Appendix B-D**).

## 2. Mathematical model

### 2.1 System dynamics

We use a common complete probability space [16]. We consider sediment storage management in a river reach. The streamflow discharge $(Q_t)_{t \geq 0}$ and lumped sediment storage $(S_t)_{t \geq 0}$ of a river reach are the state variables of our model (see **Figure 1**). Sediment storage increases only through replenishment [39].

The decision-maker, an environmental river manager, can replenish the sediment at almost surely (a.s.) strictly increasing observation times of $\tau = (\tau_k)_{k=0,1,2,...}$, where $\tau_0 = 0$, without a loss of generality. Subsequently, the sequence $\tau$ will be specified using an Erlang distribution. At each $\tau_k$ ($k \in \mathbb{N}$), the decision-maker can replenish the sediment as $S_{\tau_k+} = S_{\tau_k} + \eta_k$ with $\eta_k \geq 0$ ($k \in \mathbb{N}$), such that $S_{\tau_k} + \eta_k \leq \overline{S} < +\infty$ with a storage capacity of $\overline{S} > 0$. Here, $\eta_k = 0$ indicates no replenishment at $\tau_k$. Set the initial conditions at $S_0 \in \Omega_S \equiv [0, \overline{S}]$. Set the sequence at $\eta = (\eta_k)_{k=0,1,2,...}$ and $\eta_0 = 0$ without a loss of generality. We assume the left-continuity of $(S_t)_{t \geq 0}$ because the decision-making at $\tau_k$ should be based on the observation $(Q_{\tau_k}, S_{\tau_k})$, the consequence of which induces a new state $S_{\tau_k+}$.

Based on the relaxation technique [29], we identify the range of replenishment as $\Omega_S = \left[0, \bar{S}\right]$ and not the more natural state-dependent form $\left[0, \bar{S} - \eta_k\right]$. As the motivation for this relaxation, even if the decision-maker chooses $\eta_k$ such that $S_{\tau_k} + \eta_k > \bar{S}$, the amount of sediment effectively replenished is at most $\bar{S} - S_{\tau_k}$. This relaxation avoids using a state-dependent range while ruling out spurious optimal decisions $S_{\tau_k} + \eta_k > \bar{S}$ (Proposition 3.5 of Yoshioka et al. [29]). Thus, we have transitions $S_{\tau_k} \to \min\{S_{\tau_k} + \eta_k, \bar{S}\}$ and $Q_{\tau_k} \to Q_{\tau_k}$ (no change at $\tau_k$).

We formulate SDEs of the streamflow and sediment storage dynamics. Each SDE in each $(\tau_{k-1}, \tau_k)$ is integrated in time based on the internal conditions $(Q_{\tau_{k-1}+}, S_{\tau_{k-1}+})$ and non-negative couple $(Q_{\tau_0+}, S_{\tau_0+}) = (Q_0, S_0)$. We set a space-time Poisson random measure $N = N(\mathrm{d}z, \mathrm{d}u, \mathrm{d}t)$ on $(0, +\infty)^3$ [28] independent of $\tau$. Its compensator is $\mathrm{d}u\, v(\mathrm{d}z)\, \mathrm{d}t$ with a non-negative kernel $v = v(\mathrm{d}z)$. The compensated measure is $\tilde{N}(\mathrm{d}z, \mathrm{d}u, \mathrm{d}t) = N(\mathrm{d}z, \mathrm{d}u, \mathrm{d}t) - \mathrm{d}u\, v(\mathrm{d}z)\, \mathrm{d}t$. The discharge is assumed to follow the SDE of the gCBI type:

$$\mathrm{d}Q_t = -\rho(Q_t - \underline{Q})\mathrm{d}t + \int_0^\infty \int_0^{Q_{t-}+A} z N(\mathrm{d}z, \mathrm{d}u, \mathrm{d}t) \quad \text{for} \quad t \in (\tau_{k-1}, \tau_k) \quad (k = 1, 2, 3, \ldots) \tag{1}$$

with a constant $\rho > 0$ serving as the decay coefficient of the discharge, a (long-run) minimum discharge $\underline{Q} > 0$, and a parameter $A \geq 0$. A finite-variation tempered stable model is utilized [67], i.e.,

$$v(\mathrm{d}z) = a z^{-(1+\alpha)} e^{-bz} \mathrm{d}z \quad \text{for} \quad z > 0, \tag{2}$$

where $z$ is the jump size, $a > 0$ is the intensity parameter governing the jump intensity, $b > 0$ is the tilting parameter suppressing large jumps, and $\alpha \in [0, 1)$ is the intermittency parameter characterizing the fractional nature of small jumps. The moments of the jumps are bounded because $M_k \equiv \int_0^\infty z^k v(\mathrm{d}z) = a b^{\alpha-k} \Gamma(k-\alpha) < +\infty$ ($k \in \mathbb{N}$) with $\Gamma(\cdot)$ the Gamma function. We assume

$$\rho > M_1 \tag{3}$$

to guarantee the boundedness of the moments of the discharge as required. This is satisfied in a real case (**Section 5**). Furthermore, (3) plays a pivotal role in the proof of the comparison theorem of the HJBI equation (**Appendix A.2**).

The SDE (1) reduces to a classical CBI model if $A = 0$ [26]. The case $A > 0$ is a more flexible case in which the clustering nature of jumps can be better tuned. A classical tempered stable autoregression model [20] follows by scaling $Q_t + A \to A'Q_t + A$ with $A' \geq 0$ and then letting $A \to 0$. Any solution to (1) does not explode by Theorem 2.8(i) of Li et al. [67]. In addition, because SDE (1) is driven by an exponential decrease toward $\underline{Q}$ and positive jumps, $Q_t \geq 0$ for $t > 0$. Although this explanation is intuitive, it can be justified by Itô's formula, which shows that $\ln|Q_t|$ does not approach

$-\infty$ with a finite $t > 0$ (see Appendix A of Lungu and Øksendal [69] for a similar method). Consequently, (1) has a unique positive path-wise global solution (Theorem 3.1(i) of Li et al. [67]). The discharge domain is set at $\Omega_Q = [0, \infty)$. The natural filtration generated by $N$ up to time $t$ is denoted as $\mathcal{F}_t$. Set $\mathcal{F} = (\mathcal{F}_t)_{t \geq 0}$.

Next, we consider the sediment storage dynamics. Although the storage decreases (i.e., the sediment is flushed out from the target river reach) over time, there will be no sediment decrease if the discharge is smaller than a non-negative threshold value depending on the sediment material properties [39, 70]. We set the storage dynamics as

$$dS_t = -F(Q_t, S_t) dt \quad \text{for} \quad t \in (\tau_{k-1}, \tau_k) \quad (k = 1, 2, 3, ...), \tag{4}$$

where $F: \Omega_Q \times \Omega_S \to [0, \infty)$ is the sediment transport rate; it is non-decreasing in both arguments, $F(Q, \cdot) = 0$ for $Q \in [0, \hat{Q}]$ with a threshold $\hat{Q} > 0$, $F(Q, \cdot) > 0$ for $Q \in (\hat{Q}, \infty)$, and $F(\cdot, 0) = 0$. We require $F$ to be bounded and continuous in $\Omega_Q \times (0, \overline{S}]$. The increasing nature means that the sediment transport rate is larger for a larger flow discharge or storage [71]. The simplest model is $F(Q_t, S_t) = F_0(Q_t) \chi_{\{S_t > 0\}}$ with a non-decreasing $F_0(\cdot)$, where $\chi_{\{\cdot\}} = 1$ is the indicator function, i.e., it equals 1 when the statement inside $\{\cdot\}$ is true and is zero otherwise.

We end this subsection by introducing a natural filtration generated by discrete observations. We introduce this filtration because the decision-maker decides the interventions at each observation time $\tau_k$ using the discrete information available at this time. Following the framework of discrete and random observations [39], we set a filtration collecting discretely observed information at time $t$:

$$\mathcal{G}_t = \sigma \left\{ \left(\tau_k, Q_{\tau_k}, S_{\tau_k}\right)_{0 \leq k \leq K}, K = \sup\{k : \tau_k \leq t\} \right\}. \tag{5}$$

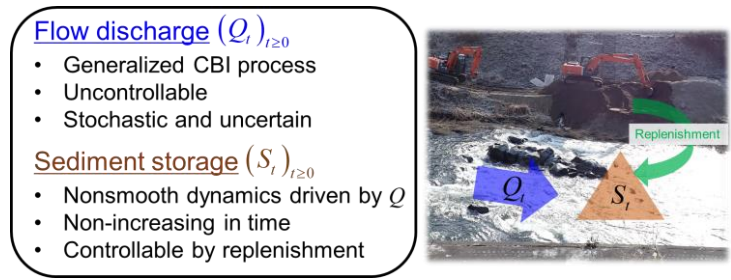

**Figure 1.** Conceptual image of the system dynamics.

## 2.2 Erlangization

Erlangization is a probabilistic technique that represents a queueing process whose waiting times follow an Erlang variable, which is semi-Markov, as a series of successive exponential variables [55]. An innovation here is that different Erlang distributions are possible at different observation times, introducing an adaptive Erlangization.

The Erlang ($l$) distribution has the probability density,

$$p_l(X) = \frac{1}{(l-1)!} \frac{X^{l-1}}{W^l} e^{-\frac{X}{W}}, \quad X \geq 0, \tag{6}$$

where $l \in \mathbb{N}$ is the shape parameter and $W > 0$ is the scaling parameter. The mean and standard deviation of the Erlang ($l$) distributions are $lW$ and $\sqrt{l}W$, respectively. The density $p_l$ is unimodal and is $X = 0$ at maximum when $l = 1$, and is $X > 0$ otherwise.

The sequence $L = (L_k)_{k=1,2,3,\ldots}$ is set with $L_k \in \Omega_L \equiv \{1, 2, 3, \ldots, \bar{L}\}$ for $\bar{L} \in \mathbb{N}$. The collection of such sequences $L = (L_k)_{k=1,2,3,\ldots}$ with each $L_k$ measurable with respect to $\mathcal{G}_{\tau_{k-1}}$ is denoted as $\mathcal{L}$. We assume that each $\Delta \tau_k = \tau_k - \tau_{k-1}$ ($k \in \mathbb{N}$) follows an Erlang ($L_k$) distribution (see **Figure 2**). For $n \in \mathbb{N}$, $\sum_{i=1}^{n} X_i \sim \text{Erlang}(n)$, where $X_1, X_2, \ldots, X_n$ are independent and identically distributed (i.i.d.) variables generated through an exponential distribution with intensity $W^{-1}$ [55]. This implies that choosing $L_k$ means choosing the total number of exponential random variables at $\tau_{k-1}$. A larger $L_k$ implies a longer waiting time.

We assume that $\eta_k$ is measurable with respect to $\mathcal{G}_{\tau_k}$, which is a natural assumption because it means that the decision-making of replenishment at $\tau_k$ is based on up-to-date information. The collection of sequences $\eta = (\eta_k)_{k=0,1,2,\ldots}$, with each $\eta_k \in \Omega_S$ ($k \in \mathbb{N}$), and $\eta_0 = 0$ for each sequence, $\eta_k$ is measurable with respect to $\mathcal{G}_{\tau_k}$, is denoted as $\mathcal{E}$. The control variables of the decision-maker are the couple $(L, \eta) \in \mathcal{L} \times \mathcal{E}$. The decision-maker therefore determines the waiting times probabilistically and can adaptively optimize the observation strategy by choosing $L_k$ at each $\tau_{k-1}$ based on $\mathcal{G}_{\tau_{k-1}}$. Our formulation naturally extends the conventional Erlangization as a special case of $L_k \equiv \bar{L}$ [55, 57–60].

*Remark 1* Here, we remark on adaptive Erlangization. The mean $\bar{X}$ is fixed, and $W = \bar{X}\bar{L}^{-1} \in \mathbb{N}$ is chosen to obtain the standard deviation $\bar{X}/\sqrt{\bar{L}}$. Letting $\bar{L} \to \infty$ implies $\text{Erlang}(\bar{L}) \to \delta_{\bar{X}}$ in the sense of distribution, where $\delta_{\bar{X}}$ is the Dirac delta at $X = \bar{X}$. The convergence speed is linear in $\bar{L}^{-1}$ measured by the coefficient of variation [72]. This result can be extended to the following: We parameterize $\bar{L} = Il' \in \mathbb{N}$ with some $I, l' \in \mathbb{N}$ and fixed values $\bar{X} = \bar{L}W$. For each $i \in \{1, 2, 3, \ldots, I\}$, the limit of $\text{Erlang}(il')$ under $l' \to \infty$ is $\delta_{\bar{X} \cdot i/I}$. The set of deterministic waiting times $\{\bar{X} \cdot i/I\}_{i=1,2,\ldots,I}$ can then be emulated through the random variables $\{X_i\}_{i=1,2,\ldots,I}$ with $X_i \sim \text{Erlang}(il')$, particularly for a large $i$. The consideration of a large $\bar{L}$ may be time-consuming in computation. Nevertheless, it provides insight into a better understanding of the relationship among the exponential, Erlang, and deterministic distributions.

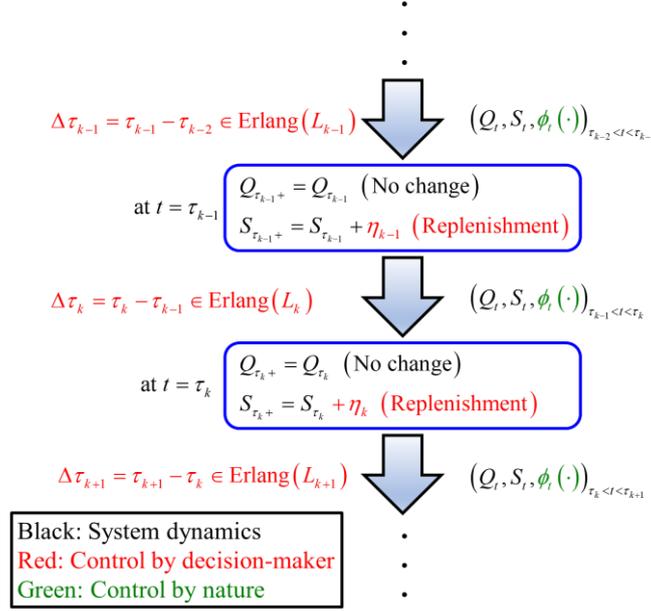

**Figure 2** Conceptual image of the adaptive Erlangization.

## 2.3 Model uncertainty

We use the convention $0\ln 0 = 0$. In our model, uncertainty is the distortion of the jump kernel, $v$, following the existing framework of robust control [61]. The probability measure of the benchmark model, i.e., the probability measure without model uncertainty, is denoted as $\mathbb{P}$. The expectation under $\mathbb{P}$ is expressed as $\mathbb{E}_{\mathbb{P}}[\cdot]$. We consider probability measures equivalent to, but possibly different from, $\mathbb{P}$.

Set a family of mapping $\phi_t(\cdot):(0,\infty)\to(0,\infty)$, $t\geq 0$, progressively measurable with respect to $\mathcal{F}$ such that

$$\mathbb{E}_{\mathbb{P}}\left[\int_0^T \int_0^\infty \int_0^{Q_s+A}(\phi_s(z)-1-\ln\phi_s(z))\mathrm{d}uv(\mathrm{d}z)\mathrm{d}s\right]<+\infty,\ T>0 \quad (7)$$

and

$$\mathbb{E}_{\mathbb{P}}\left[\exp\left\{\int_0^T \int_0^\infty \int_0^{Q_s+A}(1-\phi_s(z)+\phi_s(z)\ln\phi_s(z))\mathrm{d}uv(\mathrm{d}z)\mathrm{d}s\right\}\right]<+\infty,\ T>0. \quad (8)$$

This $\phi$ represents the model uncertainty, where $\phi_t(\cdot)\equiv 1$ means no uncertainty, and uncertainty otherwise exists. The integrands of (7) and (8) are non-negative and have the global minimum 0 at $\phi_s(z)=1$.

We set a square-integrable Martingale as a Radon–Nikodým derivative between the benchmark probability measure $\mathbb{P}$ and the distorted measure $\mathbb{Q}=\mathbb{Q}(\phi)$ parameterized by $\phi$

$$\left.\frac{\mathrm{d}\mathbb{Q}}{\mathrm{d}\mathbb{P}}\right|_t = Z_t = \exp\left\{\int_0^t \int_0^\infty \int_0^{Q_{s-}+A}\begin{pmatrix}\ln\phi_{s-}(z)\tilde{N}(\mathrm{d}z,\mathrm{d}u,\mathrm{d}s)\\-(\phi_s(z)-1-\ln\phi_s(z))\mathrm{d}uv(\mathrm{d}z)\mathrm{d}s\end{pmatrix}\right\},\ t\geq 0, \quad (9)$$

where $N(\mathrm{d}z,\mathrm{d}u,\mathrm{d}t)-\phi_t(z)\mathrm{d}u v(\mathrm{d}z)\mathrm{d}t$ becomes a compensated space-time Poisson random measure under $\mathbb{Q}$. The collection of mappings $(\phi_t(\cdot))_{t\geq 0}$ complying with (7) and (8) such that $z\min\{1,z\}\phi_t(z)v(\mathrm{d}z)$ ($z>0$) is weakly continuous, and $z^k\phi_t(z)v(\mathrm{d}z)$ ($k=1,2$) is integrable is denoted as $\mathcal{A}$. Then, $(Z_t)_{t\geq 0}$ is a martingale under $\mathbb{P}$. Weak continuity is unnecessary for a time-homogenous case ($\phi_t(\cdot)\equiv\phi(\cdot)$) (Corollary 3.2 of Kallsen and Muhle-Karbe [73]; Proposition 4.1 of Jiao et al. [66]). The set $\mathcal{A}$ is not empty because $\phi_t(z)\equiv 1$ belongs to $\mathcal{A}$. **Appendix B** provides a more generic example. Nature as the opponent of the decision-maker chooses uncertainty $\phi\in\mathcal{A}$ and thus it works against the goal of the decision-maker, as described in the next sub-section.

Finally, the relative entropy $\mathbb{D}_T(\mathbb{Q}|\mathbb{P})$ between $\mathbb{P}$ and $\mathbb{Q}$ is evaluated using the Radon–Nikodým derivative (e.g., Definition 2.1 of Hess [74]) for $T\geq 0$ as

$$\begin{aligned}\mathbb{D}_T(\mathbb{Q}|\mathbb{P}) &= \mathbb{E}_\mathbb{P}[Z_T\ln Z_T]\\ &=\mathbb{E}_\mathbb{Q}[\ln Z_T]\\ &=\mathbb{E}_\mathbb{Q}\left[\int_0^T\int_0^\infty\int_0^{Q_{s-}+A}\left(\ln\phi_{s-}(z)\tilde{N}(\mathrm{d}z,\mathrm{d}u,\mathrm{d}s)-(\phi_s(z)-1-\ln\phi_s(z))\mathrm{d}u v(\mathrm{d}z)\mathrm{d}s\right)\right].\\ &=\mathbb{E}_\mathbb{Q}\left[\int_0^T\int_0^\infty\int_0^{Q_{s-}+A}\begin{pmatrix}\ln\phi_{s-}(z)(N(\mathrm{d}z,\mathrm{d}u,\mathrm{d}s)-\phi_{s-}(z)\mathrm{d}u v(\mathrm{d}z)\mathrm{d}s)\\+(\phi_s(z)\ln\phi_s(z)-\phi_s(z)+1)\mathrm{d}u v(\mathrm{d}z)\mathrm{d}s\end{pmatrix}\right]\\ &=\mathbb{E}_\mathbb{Q}\left[\int_0^T\int_0^\infty\int_0^{Q_s+A}(\phi_s(z)\ln\phi_s(z)-\phi_s(z)+1)\mathrm{d}u v(\mathrm{d}z)\mathrm{d}s\right]\end{aligned} \quad (10)$$

## 2.4 Objective functional and effective Hamiltonian

The objective functional is an expectation as a cost-efficiency index of the observation/intervention strategy of the decision-maker subject to the worst-case uncertainty induced by nature. We consider a long-term ergodic control problem. The objective functional $J:\Omega_Q\times\Omega_S\times\mathcal{L}\times\mathcal{E}\times\mathcal{A}\to\mathbb{R}$ is set as

$$J(Q,S;L,\eta,\phi)=\limsup_{T\to+\infty}\frac{1}{T}\mathbb{E}_\mathbb{Q}\left[J_{1,T}+J_{2,T}+J_{3,T}\right], \quad (11)$$

$$J_{1,T}=\int_0^T f(S_s)\mathrm{d}s, \quad J_{2,T}=\sum_{k=1}^{\tau_k\leq T}(o+C(\eta_k)), \quad J_{3,T}=-\frac{1}{\psi}\mathbb{D}_T(\mathbb{Q}|\mathbb{P}). \quad (12)$$

The term $J_{1,T}$ measures the sediment starvation with non-negative, non-increasing, and uniformly bounded $f:\Omega_S\to[0,\infty)$ and $f\in C(0,\infty)$. This $f$ is allowed to be discontinuous at zero to penalize the depletion $S_s=0$. The term $J_{2,T}$ is the cost of the observation and intervention with the observation cost $o>0$ and replenishment cost $C:\Omega_S\to[0,\infty)$ with $C(0)=0$ (no cost for no replenishment) and $C(+0)>0$ (costly even for a small replenishment), and $C(\cdot)$ is lower semi-continuous and increasing (a larger cost for a larger amount of replenishment). A larger observation cost $o>0$ will lead to longer waiting times for optimality. The property $C(+0)>0$ implies the existence of a fixed cost [39], such as the labor costs. A standard $C$ is the piecewise linear model

$$C(\eta_k) = c_1\eta_k + c_0\chi_{\{\eta_k>0\}} \tag{13}$$

with fixed cost $c_0>0$ and the coefficient of proportional cost $c_1>0$ that we assume in what follows. The last term $J_{3,T}$ is a penalization between the benchmark and distorted models, where $\psi>0$ is the uncertainty aversion coefficient; a larger $\psi$ implies a more uncertainty-averse decision-maker [61]. Because the sum $J_{1,T}+J_{2,T}$ is the net cost to be minimized by the decision-maker, taking the minus sign in $J_{3,T}$ corresponds to a situation in which nature maximizes the net cost by paying $-J_{3,T}\geq 0$.

We set the effective Hamiltonian $H\in\mathbb{R}$ as the optimized $J$ under the uncertainty:

$$H = \inf_{(L,\eta)\in\mathcal{L}\times\mathcal{E}} \sup_{\phi\in\mathcal{A}} J(Q,S;L,\eta,\phi), \tag{14}$$

which is the saddle point of $J$. Here, we assume that the optimized objective functional is a real constant that is often called the effective Hamiltonian, which can be justified if the controlled state variables have an invariant measure. We assume its existence in this study for simplicity. Following the literature [75, 76], we justify the constant assumption (14) by verifying the HJBI equation as an optimality equation to be defined in the next subsection. A triplet of Markovian controls $(L,\eta,\phi)\in\mathcal{L}\times\mathcal{E}\times\mathcal{A}$ to give the saddle point is called an optimal control and is denoted using "*" as $(L^*,\eta^*,\phi^*)$. Owing to the Markovian nature, with the conventional abuse of notations [16], we express $L_k^* = L_k^*(Q_{\tau_{k-1}}, S_{\tau_{k-1}})$, $\eta_k^* = \eta_k^*(Q_{\tau_k}, S_{\tau_k})$, and $\phi_t^*(\cdot) = \phi^*(\cdot,Q_t,S_t)$.

**Remark 2** We see $0\leq H <+\infty$. By $J_{1,T}, J_{2,T}\geq 0$,

$$H = \inf_{(L,\eta)\in\mathcal{L}\times\mathcal{E}} \sup_{\phi\in\mathcal{A}} J(Q,S;L,\eta,\phi^1) \geq \inf_{(L,\eta)\in\mathcal{L}\times\mathcal{E}} J(Q,S;L,\eta,\phi^1) \geq 0, \tag{15}$$

where $\phi^1$ is the control with $\phi_t(z)\equiv 1$. The boundedness follows from

$$\begin{aligned}H &\leq \sup_{\phi\in\mathcal{A}} J(Q,S;L^1,\eta^0,\phi) \\ &= \sup_{\phi\in\mathcal{A}} \limsup_{T\to+\infty} \frac{1}{T}\mathbb{E}_\mathbb{Q}\left[\int_0^T f(S_s)\mathrm{d}s + \sum_{k=1}^{\tau_k\leq T} o\right] \\ &\leq \sup_{\phi\in\mathcal{A}} \limsup_{T\to+\infty} \frac{1}{T}\mathbb{E}_\mathbb{Q}\left[\int_0^T f(0)\mathrm{d}s + \sum_{k=1}^{\tau_k\leq T} o\right] \\ &= f(0) + \frac{o}{W} < +\infty \end{aligned} \tag{16}$$

where $(L^1,\eta^0)$ is the control with $L_k\equiv 1$ and $\eta_k\equiv 0$.

## 2.5 HJBI equation

Set $\Omega = \Omega_Q \times \Omega_S \times \Omega_L$. The HJBI equation is a PIDE governing the coupling $(h, \Phi)$ of $h \in \mathbb{R}$ and the potential $\Phi: \Omega \to \mathbb{R}$:

$$h + \mathbb{L}\Phi(Q,S,l) + \mathbb{M}\Phi(Q,S,l) - f(S) = 0, \quad (Q,S,l) \in \Omega. \tag{17}$$

The operator $\mathbb{L}$ acts for a generic $\Psi: \Omega \to \mathbb{R}$ as

$$\mathbb{L}\Psi(Q,S,l) = \rho(Q - \underline{Q})\frac{\partial \Psi}{\partial Q} + F(Q,S)\frac{\partial \Psi}{\partial S}$$
$$- (Q+A) \sup_{\phi(\cdot) > 0} \left\{ -\int_0^\infty \left[ \begin{matrix} (\Psi(Q,S,l) - \Psi(Q+z,S,l))\phi(z) \\ + \frac{1}{\psi}(\phi(z)\ln\phi(z) - \phi(z) + 1) \end{matrix} \right] v(\mathrm{d}z) \right\} \tag{18}$$

and the operator $\mathbb{M}$ as

$$\mathbb{M}\Psi(Q,S,l) = \frac{1}{W}\Psi(Q,S,l) - \frac{1}{W}\begin{cases} \Psi(Q,S,l-1) & (2 \le l \le L) \\ \inf_{(i,\eta) \in \Omega_L \times \Omega_S}\{\Psi(Q, \min\{S+\eta, \bar{S}\}, i) + o + C(\eta)\} & (l = 1) \end{cases}. \tag{19}$$

The operator $\mathbb{L}$ corresponds to the dynamics except at the observation times, whereas $\mathbb{M}$ corresponds to the decision-making of the observation/intervention strategy with waiting times following Erlang-type distributions chosen adaptively at each observation. The HJBI equation is satisfied over $\Omega$. In addition, the min operator appears in (19) by the relaxation technique.

The HJBI equation (19) depends on the artificial discrete parameter $l$ determining the level of Erlangization [55]. The transition from state $l$ to $l-1$ ($2 \le l \le \bar{L}$) corresponds to the jump generated by the exponential variable (See, **Section 2.2**), whereas the transition originating at $l = 1$ represents the decision-making of the Erlangization. Intuitively, the term $\mathbb{M}\Phi(\cdot, \cdot, 1)$ is related to the optimal strategy discussed below. Consequently, this HJBI equation is of a hybrid switching type with uncontrollable (the first line of (19)) and controllable (the second line) regime switching, where both the sediment replenishment and the waiting time distribution can be chosen. How the auxiliary variable $l$ effectively works can be more clearly understood in the proof of optimality (**Appendix A.1**), where it serves as an auxiliary variable of the Markovian representation of the Erlang variable.

A formal solution to the HJBI equation is a couple $(h, \Phi)$. Even if it satisfies the equation in a pointwise manner, a non-uniqueness issue arises because if $(h, \Phi)$ solves the equation, then $(h, \Phi + \Phi_0)$ with $\Phi_0 \in \mathbb{R}$ solves it as well. We add a constraint $\Phi(0,0,1) = 0$ to resolve this indeterminacy when necessary. If the HJBI equation (17) is the optimality equation and $h = H$ [75-76], then we can guess that

$$\left(L_{k+1}^*(Q_{\tau_k}, S_{\tau_k}), \eta_k^*(Q_{\tau_k}, S_{\tau_k})\right) = \arg\min_{(i,\eta) \in \Omega_L \times R(S_{\tau_k})} \{\Phi(Q_{\tau_k}, S_{\tau_k} + \eta, i) + o + C(\eta)\} \tag{20}$$

and

$$\phi^*(Q_t, S_t, \cdot) = \arg\max_{\phi(\cdot)>0} \left\{ -\int_0^\infty \begin{bmatrix} (\Phi(Q,S,\cdot) - \Phi(Q+z,S,\cdot))\phi(z) \\ +\dfrac{1}{\psi}(\phi(z)\ln\phi(z) - \phi(z) + 1) \end{bmatrix} v(\mathrm{d}z) \right\}. \quad (21)$$

$$= e^{-\psi(\Phi(Q,S,l) - \Phi(Q+z\cdot,S,l))} > 0$$

By (21), (18) becomes

$$\mathbb{L}\Phi(Q,S,l) = \rho(Q - \underline{Q})\frac{\partial \Phi}{\partial Q} + F(Q,S)\frac{\partial \Phi}{\partial S} + \frac{Q+A}{\psi}\int_0^\infty \left\{1 - e^{-\psi(\Phi(Q,S,l) - \Phi(Q+z,S,l))}\right\} v(\mathrm{d}z), \quad (22)$$

from which an exponential nonlinearity of the HJBI equation is evident.

**Remark 3** The HJBI equation (17) is redundant if the minimizing pair $(i,u) = (i^*(Q,S), u^*(Q,S))$ of $\Phi(Q, S+u, i) + o + C(u)$ satisfies $\max_{Q,S} i^*(Q,S) = \overline{i} \leq \overline{L} - 1$. In this case, the states $l = \overline{i}+1, \overline{i}+2, \ldots, \overline{L}$ are never attained, and we can omit them. However, it is difficult to determine $\overline{i}$ in advance.

## 3. Mathematical analysis

### 3.1 Viscosity solutions

Viscosity solutions as continuous weak solutions to the HJBI equation are defined based on Crandall et al. [77] and Barles and Imbert [68] because our model is a PIDE with a singular integral kernel. Hereafter, $B(Q,S,\delta)$ represents a 2-D open ball with the center $(Q,S)$ and radius $\delta > 0$, and $BUC(\Omega)$ is a collection of functions $\varphi: \Omega \to \mathbb{R}$ bounded and uniformly continuous with respect to the first and second arguments. Similarly, $C^\infty(\Omega)$ is a collection of all functions $\varphi: \Omega \to \mathbb{R}$ at differentiable arbitrary times with respect to the first and second arguments. By polynomial growth, we mean the existence of constants $\omega \in \mathbb{N}$ and $C > 0$ with $\varphi(Q, \cdot, \cdot) \leq C(1 + Q^\omega)$ for all $Q \in \Omega_Q$.

**Definition 3.1:** *A couple $(h, \Phi)$ with $h \in \mathbb{R}$ and $\Phi \in BUC(\Omega)$ is a viscosity sub-solution (resp., viscosity super-solution) if the following condition is satisfied: for any $\delta > 0$, for each $(\hat{Q}, \hat{S}, \hat{l}) \in \Omega$, and for any test function $\varphi \in C^\infty(\Omega)$ having a polynomial growth with $\Phi - \varphi$ is maximized (resp., minimized) at $(Q, S, l) = (\hat{Q}, \hat{S}, \hat{l})$ on $B(\hat{Q}, \hat{S}, \delta) \bigcap \Omega_Q \times \Omega_S$, the inequalities below are satisfied:*

$$h + \rho(\hat{Q} - \underline{Q})\frac{\partial \varphi}{\partial Q}$$
$$+ \frac{\hat{Q}+A}{\psi}\int_0^\delta \left\{1 - e^{-\psi(\varphi(\hat{Q},\hat{S},\hat{l}) - \varphi(\hat{Q}+z,\hat{S},\hat{l}))}\right\} v(\mathrm{d}z) \qquad \text{(resp., "} \geq 0 \text{") if } \hat{S} = 0 \quad (23)$$
$$+ \frac{\hat{Q}+A}{\psi}\int_\delta^\infty \left\{1 - e^{-\psi(\Phi(\hat{Q},\hat{S},\hat{l}) - \Phi(\hat{Q}+z,\hat{S},\hat{l}))}\right\} v(\mathrm{d}z) + \mathbb{M}\Phi(\hat{Q},\hat{S},\hat{l}) - f(\hat{S}) \leq 0$$

*and*

$$h + \rho(\hat{Q} - \underline{Q})\frac{\partial \varphi(\hat{Q},\hat{S},\hat{I})}{\partial Q} + F(\hat{Q},\hat{S})\frac{\partial \varphi(\hat{Q},\hat{S},\hat{I})}{\partial S}$$

$$+ \frac{\hat{Q}+A}{\psi}\int_0^\delta \left\{1 - e^{-\psi(\varphi(\hat{Q},\hat{S},\hat{I}) - \varphi(\hat{Q}+z,\hat{S},\hat{I}))}\right\}v(\mathrm{d}z) \quad \text{(resp., "} \geq 0\text{") if } \hat{S} > 0. \quad (24)$$

$$+ \frac{\hat{Q}+A}{\psi}\int_\delta^\infty \left\{1 - e^{-\psi(\Phi(\hat{Q},\hat{S},\hat{I}) - \Phi(\hat{Q}+z,\hat{S},\hat{I}))}\right\}v(\mathrm{d}z) + \mathbb{M}\Phi(\hat{Q},\hat{S},\hat{I}) - f(\hat{S}) \leq 0$$

*In addition, a couple* $(h,\Phi)$ *with* $h \in \mathbb{R}$ *and* $\Phi \in BUC(\Omega)$ *is a viscosity solution if it is a viscosity sub-solution as well as a viscosity super-solution.*

Viscosity solutions to the uncertainty-neutral case ($\psi \to +0$) are defined through a formal replacement:

$$\frac{\hat{Q}+A}{\psi}\int_0^\delta \left\{1 - e^{-\psi(\varphi(\hat{Q},\hat{S},\hat{I}) - \varphi(\hat{Q}+z,\hat{S},\hat{I}))}\right\}v(\mathrm{d}z) + \frac{\hat{Q}+A}{\psi}\int_\delta^\infty \left\{1 - e^{-\psi(\Phi(\hat{Q},\hat{S},\hat{I}) - \Phi(\hat{Q}+z,\hat{S},\hat{I}))}\right\}v(\mathrm{d}z)$$
$$\to (\hat{Q}+A)\int_0^\delta \left(\varphi(\hat{Q},\hat{S},\hat{I}) - \varphi(\hat{Q}+z,\hat{S},\hat{I})\right)v(\mathrm{d}z) + (\hat{Q}+A)\int_\delta^\infty \left(\Phi(\hat{Q},\hat{S},\hat{I}) - \Phi(\hat{Q}+z,\hat{S},\hat{I})\right)v(\mathrm{d}z) \quad (25)$$

Introducing the artificial parameter $\delta > 0$ removes the singularity of the integral of (18) at $z = 0$. Note that we assume that the equations are satisfied both along and inside the domain. This is in accord with the direction of the characteristics along the boundaries, where no information comes from outside $\Omega$.

***Remark 4*** A couple $(h,\Phi)$ with $h \in \mathbb{R}$ and $\Phi \in BUC(\Omega)$ solving the HJBI equation pointwise (classical solution) is a viscosity solution. The boundedness assumption of viscosity solutions is not restrictive because the domain is compact in the $S$ direction. Furthermore, functional (17) and hence the source term of the HJBI equation do not include $Q$. In fact, the discharge is considered to be colored noise. Removing the boundedness assumption, particularly for comparison, seems difficult because of the dependence of the integral kernel on $Q$ representing the CBI nature. Note that the boundedness assumption is standard for problems with singular jump kernels (e.g., [68]).

### 3.2 Optimality by verification

We show that the HJBI equation is the optimality equation of the control problem under regularity assumptions. The verification result (**Proposition 1**) is presented for sufficiently regular (classical) solutions, and the uniqueness result (**Proposition 2**) is presented for viscosity solutions. Therefore, the latter covers both classical and viscosity cases. Both the verification and uniqueness results are unique because of the existence of the CBI nature and model uncertainty. Note that sup in (18) and inf in (19) are attained if $\Phi \in BUC(\Omega)$.

***Proposition 1***

*Assume that there is a pair* $(h,\Phi)$ *solving the HJBI equation (17) pointwise with* $\Phi \in BUC(\Omega)$.

*Further, we assume that the couple $(Q_t, S_t)_{t \geq 0}$ admits a unique invariant measure for each $(L, \eta, \phi) \in \mathcal{L} \times \mathcal{E} \times \mathcal{A}$. We then have $h = H$. Furthermore, if (20) and (21) are admissible, they are optimal.*

### 3.3 Uniqueness

We compare the results of the HJBI equation under several assumptions. The assumption of the Hölder continuity of the sub- and super-solutions is necessary to guarantee the convergence of a non-local term representing the CBI nature. The range of the exponent of the Hölder continuity depends on $\alpha$ of the kernel $\nu$. The uncertainty aversion parameter $\psi$ must be small to manage the exponential nonlinearity of the HJBI equation. Assuming a Lipschitz continuity of the drift coefficient $F$ for $S > 0$ is rather standard, although it is discontinuous along $S = 0$. The source $f$ is allowed to be discontinuous $S = 0$. Note that even in an uncertainty-averse case, the problem in this study is non-trivial and not covered by Barles and Imbert [68] because the jump kernel is not integrable and depends on the state variable. A stronger assumption is required in the uncertainty-averse case to manage the extra nonlinearity.

*Proposition 2*

*Assume that $F(\cdot, S)$ ($S \in \Omega_S$) is Lipschitz continuous in $\Omega_Q$ and $f$ is uniformly continuous for $S > 0$.*

*(i)      Uncertainty-neutral case ($\psi = 0$)*

*For any sub-solution $(H_1, \Phi_1)$ and super-solution $(H_2, \Phi_2)$, if there is a constant $C^{(0)} > 0$ such that*

$\left| \Phi_i (Q_1, \cdot, \cdot) - \Phi_i (Q_2, \cdot, \cdot) \right| \leq C^{(0)} |Q_1 - Q_2|^\theta$ *($Q_1, Q_2 \geq 0$) with $\theta \in (\alpha, 1]$, then $H_2 \geq H_1$.*

*(ii)     Uncertainty-averse case ($\psi > 0$)*

*The statement in (i) becomes true for small $\psi$ if $\Phi_1$ and $\Phi_2$ are uniformly bounded when $\psi$ is sufficiently small.*

## 4. Numerical scheme

### 4.1 Discretization

We apply a monotone finite difference scheme to the HJBI equation. Applying monotone schemes to HJBI and related equations has been a major approach for the numerical computation of stochastic control problems in terms of viscosity solution [78]. To find $(H, \Phi)$ simultaneously, the discretized equation is solved using a fast sweeping method.

### 4.1 Computational domain

The domain $\Omega_Q$ is truncated as $[0,\bar{Q}]$ with a sufficiently large $\bar{Q}(>\hat{Q}>0)$. Similar localization techniques have been used for computing the HJBI and related equations [79, 80]. The intervals between each successive vertex in $Q$ and $S$ directions are $\Delta Q = \bar{Q}/N_Q$ and $\Delta S = \bar{S}/N_S$, respectively. Here, $N_Q+1, N_S+1 \in \mathbb{N}$ are the total number of vertices in $Q$ the and $S$ directions, respectively. The scheme aims to compute a couple $(h,\Phi)$. The former is approximated by a constant, while the latter is approximated by a collection $\Phi_{i,j,l}$ ($0 \leq i \leq N_Q$, $0 \leq j \leq N_S$, and $1 \leq l \leq \bar{L}$) of the approximated $\Phi$ at $(i\Delta Q, j\Delta S, l)$.

## 4.2 Discretization

We first describe the discretization for an uncertainty-neutral case followed by that for an uncertainty-averse case. First, the integral term is rewritten as follows:

$$\int_0^\infty \left(\Phi(Q,S,l)-\Phi(Q+z,S,l)\right)v(\mathrm{d}z) = \int_0^\infty \left(\Phi(Q,S,l)-\Phi(Q+z,S,l)+z\frac{\partial \Phi}{\partial Q}\right)v(\mathrm{d}z) - M_1\frac{\partial \Phi}{\partial Q}, \quad (26)$$

through which the HJBI equation in an uncertainty-neutral case becomes

$$\begin{aligned}
& h - f(S) \\
& + \left((\rho-M_1)Q - AM_1 - \rho\underline{Q}\right)\frac{\partial \Phi}{\partial Q} + F(Q,S)\frac{\partial \Phi}{\partial S} \\
& + (Q+A)\int_0^\infty \left(\Phi(Q,S,l)-\Phi(Q+z,S,l)+z\frac{\partial \Phi}{\partial Q}\right)v(\mathrm{d}z) \\
& + \mathbb{M}\Phi = 0
\end{aligned} \quad (27)$$

Thus, the integrand of the nonlocal term becomes formally smoother near $z=0$.

Each line in (27) is discretized as follows: The first line is discretized as $h-f(S_j)$. The second line is based on a one-sided finite difference:

$$\left((\rho-M_1)Q - AM_1 - \rho\underline{Q}\right)\frac{\partial \Phi}{\partial Q}$$

$$\rightarrow \left((\rho-M_1)Q_i - AM_1 - \rho\underline{Q}\right) \times \begin{cases} \dfrac{\Phi_{i,j,l}-\Phi_{i-1,j,l}}{\Delta Q} & \left((\rho-M_1)Q_i - AM_1 - \rho\underline{Q} \geq 0,\ i < N_Q\right) \\ \dfrac{\Phi_{i+1,j,l}-\Phi_{i,j,l}}{\Delta Q} & \left((\rho-M_1)Q_i - AM_1 - \rho\underline{Q} < 0,\ i < N_Q\right) \\ 0 & (i = N_Q) \end{cases}, \quad (28)$$

$$F(Q,S)\frac{\partial \Phi}{\partial S} \rightarrow \begin{cases} F(Q_i,S_j)\dfrac{\Phi_{i,j,l}-\Phi_{i,j-1,l}}{\Delta S} & (j>0) \\ 0 & (j=0) \end{cases}. \quad (29)$$

The third line is discretized as

$$(Q+A)\int_0^\infty \left[\Phi(Q,S,l)-\Phi(Q+z,S,l)+z\frac{\partial\Phi}{\partial Q}\right]v(\mathrm{d}z)$$

$$\rightarrow \begin{cases} (Q_i+A)\sum_{k=1}^{N_Q-i} I_{i,j,l,k} + (A+Q_i)V\dfrac{\Phi_{i+1,j,l}-\Phi_{i,j,l}}{\Delta Q} & (i=0) \\ (Q_i+A)\sum_{k=1}^{N_Q-i} I_{i,j,l,k} & (0<i<N_Q) \\ 0 & (i=N_Q) \end{cases} \quad (30)$$

with

$$I_{i,j,l,k}=v_k\left(\Phi_{i,j,l}-\Phi_{i+k,j,l}+\chi_{\{i>0\}}z_k\frac{\Phi_{i,j,l}-\Phi_{i-1,j,l}}{\Delta Q}\right),\quad v_k=a\frac{\exp(-b(k+1/2)\Delta Q)}{((k+1/2)\Delta Q)^{1+\alpha}}\Delta Q, \quad (31)$$

$$V=\int_{\bar{Q}}^\infty zv(\mathrm{d}z)=a\bar{Q}^{1-\alpha}\int_0^1 x^{\alpha-2}e^{-\frac{b\bar{Q}}{x}}\mathrm{d}x \in (0,M_1), \quad (32)$$

where $\sum_{k=1}^{N_Q-i} I_{i,j,l,k}$ is the discretized integral term, and $v_k$ is the integration of $v$. The integral (32) is computed using the classical Simpson rule. Note that $\Phi_{i,j,l}-\Phi_{i-1,j,l}$ cannot be defined for $i=0$.

Finally, the fourth line is discretized as

$$\mathbb{M}\Phi \rightarrow \frac{1}{W}\Phi_{i,j,l}-\frac{1}{W}\begin{cases}\Phi_{i,j,l-1} & (2\le l\le L) \\ \min_{(j',l')}\left\{\Phi_{i,\min\{j+j',N_S\},l'}+o+C(j'\Delta S)\right\} & (l=1)\end{cases}, \quad (33)$$

where the minimization with respect to non-negative integers $(j',l')$ is applied for $j'\in\{0,1,2,...,N_S\}$ and $l'\in\{1,2,...,\bar{L}\}$. Substituting (28)–(33) into (27) yields the nonlinear system of the form

$$h-f(S_j)+\Xi_{i,j,l}\Phi_{i,j,l}-\Theta_{i,j,l}[\Phi]=0 \text{ for } 0\le i\le N_Q,\ 0\le j\le N_S,\ 1\le l\le \bar{L}, \quad (34)$$

where $\Xi_{i,j,l}>0$, $\Theta_{i,j,l}[\Phi]$ is a summation of the terms containing $\Phi_{i+i',j+j',l+l'}$ with $i',j',l'$ such that $i'j'l'\ne 0$. Based on the construction, $\Theta_{i,j,l}[\Phi]$ is nondecreasing with respect to each $\Phi_{i+i',j+j',l+l'}$. We have $(N_Q+1)(N_S+1)\bar{L}+1$ unknowns but only $(N_Q+1)(N_S+1)\bar{L}$ equations, and the indeterminacy is met. Thus, we chose $\Phi_{0,0,0}=0$ to resolve this issue. The discretization for $i=0$ of (28) and (30) satisfies the monotonicity as a whole because

$$\left((\rho-M_1)\cdot 0 - AM_1-\rho\underline{Q}\right)\frac{\Phi_{0+1,j,l}-\Phi_{0,j,l}}{\Delta Q}+(0+A)\sum_{k=1}^{N_Q-0} I_{0,j,l,k}+(A+0)V\frac{\Phi_{0+1,j,l}-\Phi_{0,j,l}}{\Delta Q}$$
$$=A\sum_{k=1}^{N_Q} I_{0,j,l,k}+\left(A(M_1-V)+\rho\underline{Q}\right)\frac{\Phi_{0,j,l}-\Phi_{1,j,l}}{\Delta Q}. \quad (35)$$

A similar argument is true for the uncertainty-averse case discussed below. The condition $\Xi_{i,j,l}>0$ justifies the fast-sweeping method described in the next subsection.

For an uncertainty-averse case, the non-local term is discretized as

$$\frac{Q+A}{\psi}\int_0^\infty \left\{1-e^{-\psi(\Phi(Q,S,l)-\Phi(Q+z,S,l))}\right\}v(\mathrm{d}z)$$

$$=\frac{Q+A}{\psi}\int_0^\infty \left\{1-e^{-\psi(\Phi(Q,S,l)-\Phi(Q+z,S,l))}+\psi z\frac{\partial \Phi}{\partial Q}\right\}v(\mathrm{d}z)-M_1(Q+A)\frac{\partial \Phi}{\partial Q}$$

$$\approx \frac{Q+A}{\psi}\int_0^{\bar{Q}} \left\{1-e^{-\psi(\Phi(Q,S,l)-\Phi(Q+z,S,l))}+\psi z\frac{\partial \Phi}{\partial Q}\right\}v(\mathrm{d}z)-M_1(Q+A)\frac{\partial \Phi}{\partial Q}$$

$$\to \begin{cases} \dfrac{Q_i+A}{\psi}\sum_{k=1}^{N_Q-i}J_{i,j,l,k}+(A+Q_i)V\dfrac{\Phi_{i+1,j,l}-\Phi_{i,j,l}}{\Delta Q} & (i=0) \\ \dfrac{Q_i+A}{\psi}\sum_{k=1}^{N_Q-i}J_{i,j,l,k} & (0<i<N_Q) \\ 0 & (i=N_Q) \end{cases} \quad (36)$$

with

$$J_{i,j,l,k}=v_k\left(1-e^{-\psi(\Phi_{i,j,l}-\Phi_{i+k,j,l})}+\psi\chi_{\{i>0\}}z_k\frac{\Phi_{i,j,l}-\Phi_{i-1,j,l}}{\Delta Q}\right). \quad (37)$$

This $J_{i,j,l,k}$ is increasing with respect to $\Phi_{i,j,l}-\Phi_{i+k,j,l}$ and $\Phi_{i,j,l}-\Phi_{i-1,j,l}$, serving as a key in the boundedness and uniqueness of numerical solutions (**Appendices A.3-A.4**). The discretized HJBI equation has the form

$$h-f(S_j)+\Xi_{i,j,l}\Phi_{i,j,l}-\Theta_{i,j,l}[\Phi]+\Lambda_{i,j,l}[\Phi]=0 \quad (38)$$

with

$$\Lambda_{i,j,l}[\Phi]=\chi_{\{i<N_Q\}}\frac{Q_i+A}{\psi}\sum_{k=1}^{N_Q-i}v_k\left(1-e^{-\psi(\Phi_{i,j,l}-\Phi_{i+k,j,l})}-\psi(\Phi_{i,j,l}-\Phi_{i+k,j,l})\right). \quad (39)$$

Finally, discretized optimal controls at each $(i\Delta Q, j\Delta S, l)$ are found as

$$\left(\frac{\eta^*(Q_i,S_j)}{\Delta S},L^*(Q_i,S_j)\right)=\arg\min_{(j',l')}\left\{\Phi_{i,\min\{j+j',N_S\},l'}+o+C(j'\Delta S)\right\}, \quad (40)$$

$$\phi^*(Q_i,S_j,i'\Delta Q)=\exp\left(-\psi(\Phi_{i,j,l}-\Phi_{i+i',j,l})\right) \quad (i'=0,1,2,...,N_Q-i). \quad (41)$$

*Remark 5* Higher-order discretization methods, such as the (weighted) essentially non-oscillatory schemes, can be applied to first-order derivatives [81, 82]. However, we did not find significant improvement. This is considered to be due to the dominant errors induced by the non-local terms. In addition, **Propositions 3** and **4** below seem not to hold true for such schemes.

The following propositions provide insight into numerical solutions. **Proposition 3** shows their bound complying with the continuous counterpart in **Remark 2**, and **Proposition 4** shows the uniqueness of numerical solutions. As the novelty of **Proposition 4**, our discretization is not of a proper type by which the uniqueness of the solutions follows under mild conditions (Definition 6 of Oberman [83]). We show that monotone discretization combined with the connectivity assumption of the scheme overcomes this difficulty. Hereafter, the range of subscripts $(i,j,l)$ are omitted when there will be no confusion.

*Proposition 3*

*If there is a numerical solution $\left(h,\{\Phi_{i,j,l}\}\right)$ to (38) or (34) with each $\Phi_{i,j,l}$ bounded, then*

$$0 \leq h \leq f(0) + oW^{-1}. \tag{42}$$

*Proposition 4*

*Assume $\Phi_{0,0,0} = \Psi_{0,0,0} = 0$. For sufficiently small $\bar{L}$ and $\bar{S}$, for any two numerical solutions $\left(\bar{h},\{\Phi_{i,j,l}\}\right)$ and $\left(\underline{h},\{\Psi_{i,j,l}\}\right)$ of (38) or (34) with $\bar{h} \geq \underline{h}$, we have $\{\Phi_{i,j,l}\} \leq \{\Psi_{i,j,l}\}$ component wise. Furthermore, if $\bar{h} = \underline{h}$, then $\{\Phi_{i,j,l}\} \equiv \{\Psi_{i,j,l}\}$.*

**Remark 6** It seems difficult to bound $\Phi_{i,j,l}$ in an *a priori* manner as in the discounted cases (Lemma 5.1 of Forsyth and Labahn [84]) because our system is not strictly monotone. Monotone and consistent schemes converge in a viscosity sense for discounted cases [78]. Our scheme is neither proper nor Lipschitz (Definition 3 of Oberman [83]) in an uncertainty-averse case. If one replaces $h$ by $\delta \Phi_{i,j,l}$, with a small $\delta > 0$ as a small-discount counterpart, the unique existence of the numerical solutions will be proven; however, the problem is now perturbed. The lack of the aforementioned properties hinders the application of conventional existence arguments for numerical solutions. We checked the performance of the proposed scheme against text cases before applying it to realistic cases.

**Remark 7** For an uncertainty-neutral case, assuming the existence of numerical solutions is not restrictive. This is because our scheme is monotone and all vertices are connected (see **Proof of Proposition 4**); thus, the discretized HJBI equation is seen as a Bellman equation of a discrete-time ergodic Markov chain. We then apply Theorem 7.58 in [85] to determine whether the chain has an invariant measure. Moreover, the existence of a numerical solution implies an optimality of the controls (40) in a discrete sense (Theorem 6.17 in [86]) owing to their countable nature.

### 4.3 Fast-sweeping method

We solve the discretized equation using the fast-sweeping method of Yoshioka et al. [87]. This method can avoid inverting the coefficient matrix, which is dense in our case. The numerical solution at the $n$ th sweeping step is represented by the super-script $(n)$. Starting from the initial guess $\left(h^{(0)},\{\Phi_{i,j,l}^{(0)}\}\right) \equiv (0,0)$, the iteration proceeds as follows. The quantity $\varepsilon > 0$ is the error tolerance used to judge the convergence of the iteration. The weight $w \in (0,1)$ is introduced for a stable computation.

**Step 1.** Compute $h^{(n+1)} = f(S_j) - \Theta_{i,j,l}\left[\Phi^{(n)}\right] - \Lambda_{i,j,l}\left[\Phi^{(n)}\right] = 0$ and go to **Step 2**.

**Step 2.** Compute $\Phi_{i,j,l}^{(n+1)} = (1-w)\Phi_{i,j,l}^{(n)} + w\frac{1}{\Xi_{i,j,l}}\left\{f(S_j) - h^{(n+1)} + \Theta_{i,j,l}\left[\Phi^{(n)}\right] + \Lambda_{i,j,l}\left[\Phi^{(n)}\right]\right\}$ for all $i,j,l$ except at $(i,j,l) = (0,0,0)$ and go to **Step 3**.

**Step 3.** The iteration error is computed as $\mathrm{Er} = \max_{i,j,l}\left|\Phi_{i,j,l}^{(n+1)} - \Phi_{i,j,l}^{(n)}\right|$. If $\mathrm{Er} > \varepsilon$, then $n \to n+1$ and proceed to **Step 1**. Otherwise, the output $\left(h^{(n+1)}, \left\{\Phi_{i,j,l}^{(n+1)}\right\}\right)$ is used as the numerical solution.

Here, $\Phi_{i,j,l}^{(n)}$ is the latest component to be used for calculating the terms $\Theta_{i,j,l}\left[\Phi^{(n)}\right]$, $\Lambda_{i,j,l}\left[\Phi^{(n)}\right]$ following the Gauss–Seidel sweeping [88]. We set $w = 0.35$ for test cases, and $w = 0.30$ otherwise.

## 5 Application
### 5.1 Parameter identification
#### 5.1.1 Streamflow

We specify the parameter values of the SDEs. We used the public hourly hydrological data of Obara Dam, the largest dam along Hii River in Japan [89]. Its downstream river reach has been a study site for sediment replenishment strategies, and a sediment replenishment project with $\bar{S} = O(10^2)$ (m³) was initiated in April 2020 [39]. Its optimization is currently an urgent issue for this river. We utilized the outflow discharge data from April 1, 2016, to March 31, 2020 (**Figure 3**) to identify the streamflow dynamics. The minimum discharge corresponding to the maintenance flow is 1.0 m³/s, suggesting that $\underline{Q} = 1.0$ (m³/s) should be set. The dam discharge is clustered over time, as shown in **Figure 3**.

The streamflow dynamics have five parameters: $a, b, \alpha, A, \rho$. Set $a' = a\rho$. First, we specify the four parameters $a', b, \alpha, A$ using a method of moment matching similar to Yoshioka and Yoshioka [20] such that the following error metric is minimized: $\mathrm{Er} = \sum_{i=1}^{4}\left((c_{i,\mathrm{O}} - c_{i,\mathrm{M}})/c_{i,\mathrm{O}}\right)^2$, where $c_1, c_2, c_3, c_4$ represents the stationary average, standard deviation, skewness, and kurtosis, respectively (**Appendix D**). The subscripts O and M represent the observed and modeled quantities, respectively. Using a common nonlinear least-squares method, we obtain $\alpha = 0.201$ (-), $a' = 3.49\times 10^{-3}$ ( $\mathrm{m}^{3(\alpha-1)}/\mathrm{s}^{\alpha-1}$ ), $b = 8.33\times 10^{-3}$ (s/m³), $A = 16.5$ (m³/s), and the moments in **Table 1** with $\mathrm{Er} = 4.52\times 10^{-9}$. The relative errors of the average and standard deviations are less than 0.005%, suggesting an excellent agreement between the observations and the model. From the identification results, we obtain $M_1 = 0.187\rho$. For the remaining parameter $\rho$, we use the autocorrelation function of the gCBI process, which in our case is an exponential function $\exp(-(\rho - M_1)\delta)$ with a time lag $\delta \geq 0$ (**Appendix C**). We identified the correlation coefficient $\rho - M_1 = 0.813\rho = 0.028$ (1/h) and hence $\rho = 0.033$ (1/h).

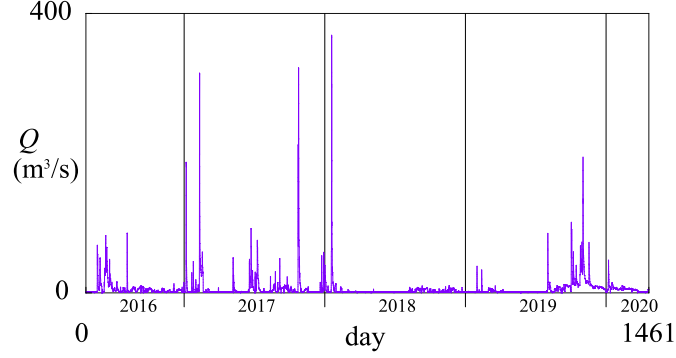

**Figure 3.** Time series data of the outflow discharge from April 1, 2016 to March 31, 2020.

**Table 1.** Comparison of the observed and modelled statistics of the discharge.

| Moments | Observed | Modeled | Relative error |
|---|---|---|---|
| Average | 5.014.E+00 | 5.014.E+00 | 3.654.E-05 |
| Standard deviation | 1.539.E+01 | 1.539.E+01 | 3.066.E-05 |
| Skewness | 1.198.E+01 | 1.198.E+01 | 3.799.E-05 |
| Kurtosis | 1.980.E+02 | 1.980.E+02 | 2.834.E-05 |

### 5.1.2 Sediment

For the coefficient $F$, we use the physically based sediment transport formula, which has been extensively verified against a wide range of flow conditions [90]. The formula calculates the sediment transport rate $Q^{(s)}$ and momentum (mass times velocity) of sediment per unit of riverbed area as

$$\frac{Q^{(s)}}{\rho_p d \sqrt{\left(\frac{\rho_p}{\rho_w} - 1\right)gd}} = \frac{2\sqrt{\theta_t}}{\kappa \mu_b \sqrt{Z_c}} \max\{\theta - \theta_t, 0\}\left[1 + \frac{c_M}{\mu_b}\max\{\theta - \theta_t, 0\}\right] \text{ (kg m/s/m}^2\text{).} \quad (43)$$

Here, $\rho_p$ is the particle density (2,650 kg/m³), $\rho_w$ is the water density (997 kg/m³), $g$ is the gravitational acceleration (9.81 m/s²), $d$ is the sediment diameter in m, $\kappa$ is the Karman constant (0.4-), $\mu_b$ is the empirical coefficient (0.63 (-)), $c_M$ is another coefficient (1.7 (-)), $\theta$ is the Shields number (non-dimensional bed shear stress (-)), $\theta_t$ is the threshold Shields number (approximately 0.07 (-)), and $Z_c = 1 - \mu_b^{-1} I_w \rho_p / (\rho_p - \rho_w)$ is the slope correction (-). The Shields number $\theta$ is calculated from the bed shear stress $\tau_w = \rho_w g H_w I_w$ (where $H_w$ is the water depth (m), and $I_w \ll 1$ is the bed slope (-)) as follows:

$$\theta = \frac{\tau_w}{(\rho_p - \rho_w)gdZ_c}. \quad (44)$$

Furthermore, for a wide rectangular channel cross-section, the water depth $H_w$ as a function of the discharge $Q$ can be estimated using the Manning formula [91]:

$$H_w = H_w(Q) = \left(n_w Q B_w^{-1} I_w^{-0.5}\right)^{0.6} \text{ (m)}, \tag{45}$$

where $B_w$ is the channel width (m) and $n_w$ is the roughness coefficient. We then obtain the sediment transport rate $Q^{(s)}$ as an increasing function of the discharge $Q$ based on (43)–(45). Assume that the sediment storage is measured in units of volume, and its porosity is $\varsigma \in (0,1)$. Then, $F$ is related to $Q^{(s)}$ as

$$F = \chi_{\{S>0\}} \frac{B_w Q^{(s)}}{\varsigma \rho_p} \text{ (m}^3\text{/s)}. \tag{46}$$

In summary, with the truncation for a numerical computation, $F$ is a function of $Q$ in the form

$$F(Q,S) = \chi_{\{S>0\}} F_1 \max\left\{\min\{Q,\bar{Q}\}^{0.6} - \hat{Q}^{0.6}, 0\right\}\left(1 + F_2 \max\left\{\min\{Q,\bar{Q}\}^{0.6} - \hat{Q}^{0.6}, 0\right\}\right) \tag{47}$$

with $F_0, F_1, \hat{Q} > 0$. This $F$ is discontinuous in $S$ and is Hölder continuous with an exponent of 0.6 in $Q$. This coefficient does not satisfy the assumption of the comparison theorem and is therefore not covered by the theory. Nevertheless, the computation below suggests the convergence of the numerical solutions.

## 5.2 Test case with manufactured solutions

The computational performance of the proposed numerical scheme is examined against manufacture solutions to a reduced model. Because the potential difficulty in the numerical computation of the HJBI equation will be a non-local term having a singular jump kernel, we consider the following equation:

$$\begin{aligned} & H + \rho(Q - \underline{Q})\frac{\partial \Phi}{\partial Q} + F(Q,S)\frac{\partial \Phi}{\partial S} \\ & + (Q+A)\int_0^\infty (\Phi(Q,S) - \Phi(Q+z,S))\nu(\mathrm{d}z) - \bar{f}(Q,S) = 0 \end{aligned}, \quad (Q,S) \in \Omega_Q \times \Omega_S \tag{48}$$

with a bivariate function $\bar{f} = \bar{f}(Q,S)$. We choose $\bar{f}$ such that (48) admits the following solution:

$$H = 1 \text{ and } \Phi(Q,S) = -AQS^\beta + \Phi_0 \tag{49}$$

with constants $A, \Phi_0 \in \mathbb{R}$ and $\beta > 0$. The couple $(H, \Phi)$ is a viscosity solution to (48) without boundedness. Note that in the original HJBI equation, $f$ does not depend on $Q$. We set $\Phi(0,0) = 0$. This solution is simple but serves as an efficient case to check the performance of the numerical scheme because its regularity, which is expected to affect the convergence rate of numerical solutions, can be tuned by $\beta$. The test case here is interesting because it is a degenerate elliptic problem with a discontinuous source term having a discontinuous drift and source terms along the boundary $S = 0$. We choose $N_Q = N_S$ and use normalization $(Q, S) \to \left(\dfrac{Q}{\bar{Q}}, \dfrac{S}{\bar{S}}\right)$ to obtain the normalized domain $[0,1]^2$.

We use the following parameter values in dimensional form for the manufactured solutions: $A = 1$ (−), $\Phi_0 = 0$ (−), $d = 0.005$ (m), $\varsigma = 0.5$ (−), $n_w = 0.035$ (s/m$^{1/3}$), $B_w = 20$ (m), $I_w = 0.0015$ (-), $\theta_t = 0.072$ (-), and $\bar{S} = 400$ (m$^3$). In addition, $\bar{f}$ is discretized as $\bar{f}(Q_i, S_j)$. We choose a sufficiently large $\bar{Q} = 200$ (m$^3$/s) such that the probability of occurrence of the events $Q > \bar{Q}$ is less than 0.1%. We examine three cases $\beta = 0.5$ (concave case), $\beta = 1.0$ (linear case), and $\beta = 2.0$ (convex case) and use a small tolerance of $\varepsilon = 10^{-10}$ by which the sweeping is terminated at $O(10^3)$ iterations.

**Tables 2–4** show the computational errors ($l^1$, $l^2$, and $l^\infty$ norms for $\Phi$ and the absolute error for $H$) and the corresponding convergence rates, suggesting the convergence of the scheme. A faster convergence is achieved for larger $\beta$ values with a better regularity of the manufactured solution. The presented scheme is experimentally first-order accurate for a convex case and almost 0.5th order accurate for a concave case. These convergence rates are not critically different from the existing results (e.g., [92, 93]). A higher convergence rate is impossible because of the use of monotone discretization. The convergence rate is smaller for the case with a less regular solution, as expected. The convergence rate of $H$ seems to be small in some cases; however, its computational error is already small for $N_Q = 10$. The errors of the linear case are dominated by a domain truncation, as suggested in **Tables 5–7** for $\bar{Q} = 300$ (m$^3$/s). The convergence rates in the convex and concave cases were unaffected by the truncation.

**Table 2.** Computational error and convergence rate for a concave case ($\beta = 0.5$, $\bar{Q} = 200$).

| $N_Q$ | Error | | | | Convergence rate | | | |
|---|---|---|---|---|---|---|---|---|
| | $H$ | $l^1(\Phi)$ | $l^2(\Phi)$ | $l^\infty(\Phi)$ | $H$ | $l^1(\Phi)$ | $l^2(\Phi)$ | $l^\infty(\Phi)$ |
| 10 | -1.34.E-05 | 3.95.E-02 | 5.03.E-02 | 9.04.E-02 | 9.26.E-01 | 4.04.E-01 | 4.61.E-01 | 4.71.E-01 |
| 20 | -7.03.E-06 | 2.99.E-02 | 3.65.E-02 | 6.53.E-02 | 1.66.E+00 | 4.37.E-01 | 4.66.E-01 | 4.72.E-01 |
| 40 | -2.22.E-06 | 2.20.E-02 | 2.64.E-02 | 4.71.E-02 | 1.24.E+00 | 4.58.E-01 | 4.75.E-01 | 4.79.E-01 |
| 80 | -9.41.E-07 | 1.61.E-02 | 1.90.E-02 | 3.38.E-02 | 1.36.E+00 | 4.71.E-01 | 4.81.E-01 | 4.84.E-01 |
| 160 | -3.66.E-07 | 1.16.E-02 | 1.36.E-02 | 2.41.E-02 | 6.94.E-01 | 4.78.E-01 | 4.85.E-01 | 4.87.E-01 |
| 320 | -2.26.E-07 | 8.32.E-03 | 9.74.E-03 | 1.72.E-02 | | | | |

**Table 3.** Computational error and convergence rate for a linear case ($\beta = 1.0$, $\bar{Q} = 200$).

| $N_Q$ | Error | | | | Convergence rate | | | |
|---|---|---|---|---|---|---|---|---|
| | $H$ | $l^1(\Phi)$ | $l^2(\Phi)$ | $l^\infty(\Phi)$ | $H$ | $l^1(\Phi)$ | $l^2(\Phi)$ | $l^\infty(\Phi)$ |
| 10 | 8.34.E-07 | 3.75.E-04 | 5.06.E-04 | 1.42.E-03 | -7.00.E-01 | 2.93.E-01 | 3.00.E-01 | -1.92.E-01 |
| 20 | 1.36.E-06 | 3.06.E-04 | 4.11.E-04 | 1.62.E-03 | -1.26.E+00 | 4.67.E-01 | -1.27.E-01 | -8.93.E-01 |
| 40 | 3.24.E-06 | 2.22.E-04 | 4.48.E-04 | 3.01.E-03 | 6.26.E-02 | 4.71.E-01 | 3.33.E-01 | 1.61.E-01 |
| 80 | 3.10.E-06 | 1.60.E-04 | 3.56.E-04 | 2.69.E-03 | 6.59.E-02 | 2.75.E-01 | 2.14.E-01 | 1.62.E-01 |
| 160 | 2.96.E-06 | 1.32.E-04 | 3.07.E-04 | 2.40.E-03 | 9.62.E-02 | 8.93.E-02 | 2.19.E-01 | 2.18.E-01 |
| 320 | 2.77.E-06 | 1.24.E-04 | 2.64.E-04 | 2.07.E-03 | | | | |

**Table 4.** Computational error and convergence rate for a convex case ($\beta = 2.0$, $\bar{Q} = 200$).

| $N_Q$ | Error | | | | Convergence rate | | | |
|---|---|---|---|---|---|---|---|---|
| | $H$ | $l^1(\Phi)$ | $l^2(\Phi)$ | $l^\infty(\Phi)$ | $H$ | $l^1(\Phi)$ | $l^2(\Phi)$ | $l^\infty(\Phi)$ |
| 10 | 1.83.E-07 | 7.36.E-03 | 1.02.E-02 | 2.90.E-02 | -3.96.E+00 | 1.00.E+00 | 1.06.E+00 | 1.04.E+00 |
| 20 | 2.85.E-06 | 3.68.E-03 | 4.89.E-03 | 1.41.E-02 | -3.55.E-01 | 1.01.E+00 | 1.05.E+00 | 1.04.E+00 |
| 40 | 3.64.E-06 | 1.83.E-03 | 2.36.E-03 | 6.82.E-03 | 2.05.E-01 | 1.03.E+00 | 1.04.E+00 | 1.04.E+00 |
| 80 | 3.16.E-06 | 8.94.E-04 | 1.15.E-03 | 3.32.E-03 | 1.00.E-01 | 1.09.E+00 | 1.04.E+00 | 5.17.E-01 |
| 160 | 2.95.E-06 | 4.21.E-04 | 5.56.E-04 | 2.32.E-03 | 1.11.E-01 | 1.18.E+00 | 9.97.E-01 | 1.95.E-01 |
| 320 | 2.73.E-06 | 1.86.E-04 | 2.79.E-04 | 2.02.E-03 | | | | |

**Table 5.** Computational error and convergence rate for a concave case ($\beta = 0.5$, $\bar{Q} = 300$).

| $N_Q$ | Error | | | | Convergence rate | | | |
|---|---|---|---|---|---|---|---|---|
| | $H$ | $l^1(\Phi)$ | $l^2(\Phi)$ | $l^\infty(\Phi)$ | $H$ | $l^1(\Phi)$ | $l^2(\Phi)$ | $l^\infty(\Phi)$ |
| 10 | -1.07.E-05 | 4.02.E-02 | 5.10.E-02 | 9.19.E-02 | 8.40.E-01 | 4.03.E-01 | 4.63.E-01 | 4.76.E-01 |
| 20 | -5.96.E-06 | 3.04.E-02 | 3.70.E-02 | 6.61.E-02 | 1.15.E+00 | 4.38.E-01 | 4.69.E-01 | 4.78.E-01 |
| 40 | -2.68.E-06 | 2.24.E-02 | 2.67.E-02 | 4.75.E-02 | 1.53.E+00 | 4.58.E-01 | 4.74.E-01 | 4.81.E-01 |
| 80 | -9.25.E-07 | 1.63.E-02 | 1.92.E-02 | 3.40.E-02 | 1.80.E+00 | 4.73.E-01 | 4.82.E-01 | 4.85.E-01 |
| 160 | -2.66.E-07 | 1.18.E-02 | 1.38.E-02 | 2.43.E-02 | 2.33.E+00 | 4.81.E-01 | 4.86.E-01 | 4.89.E-01 |
| 320 | -5.27.E-08 | 8.42.E-03 | 9.84.E-03 | 1.73.E-02 | | | | |

**Table 6.** Computational error and convergence rate for a linear case ($\beta = 1.0$, $\bar{Q} = 300$).

| $N_Q$ | Error | | | | Convergence rate | | | |
|---|---|---|---|---|---|---|---|---|
| | $H$ | $l^1(\Phi)$ | $l^2(\Phi)$ | $l^\infty(\Phi)$ | $H$ | $l^1(\Phi)$ | $l^2(\Phi)$ | $l^\infty(\Phi)$ |
| 10 | 4.49.E-07 | 2.77.E-04 | 3.69.E-04 | 9.60.E-04 | -5.61.E-01 | 3.04.E-01 | 3.64.E-01 | 3.50.E-03 |
| 20 | 6.62.E-07 | 2.24.E-04 | 2.87.E-04 | 9.58.E-04 | -8.37.E-01 | 5.09.E-01 | 3.51.E-01 | -4.01.E-01 |
| 40 | 1.18.E-06 | 1.58.E-04 | 2.25.E-04 | 1.27.E-03 | -5.71.E-01 | 7.32.E-01 | 2.17.E-01 | -2.01.E-01 |
| 80 | 1.76.E-06 | 9.49.E-05 | 1.94.E-04 | 1.45.E-03 | 5.13.E-03 | 3.52.E-01 | 2.08.E-01 | -4.11.E-02 |
| 160 | 1.75.E-06 | 7.44.E-05 | 1.68.E-04 | 1.50.E-03 | 1.18.E-01 | 1.77.E-01 | 1.72.E-01 | 1.40.E-01 |
| 320 | 1.61.E-06 | 6.58.E-05 | 1.49.E-04 | 1.36.E-03 | | | | |

**Table 7.** Computational error and convergence rate for a convex case ($\beta = 2.0$, $\bar{Q} = 300$).

| | | Error | | | | Convergence rate | | |
|---|---|---|---|---|---|---|---|---|
| $N_Q$ | $H$ | $l^1(\Phi)$ | $l^2(\Phi)$ | $l^\infty(\Phi)$ | $H$ | $l^1(\Phi)$ | $l^2(\Phi)$ | $l^\infty(\Phi)$ |
| 10 | 3.86.E-06 | 7.40.E-03 | 1.02.E-02 | 2.91.E-02 | 1.01.E+00 | 1.01.E+00 | 1.06.E+00 | 1.05.E+00 |
| 20 | 1.92.E-06 | 3.68.E-03 | 4.90.E-03 | 1.41.E-02 | 3.93.E-01 | 9.95.E-01 | 1.03.E+00 | 1.02.E+00 |
| 40 | 1.46.E-06 | 1.85.E-03 | 2.40.E-03 | 6.93.E-03 | -3.12.E-01 | 1.05.E+00 | 1.06.E+00 | 1.04.E+00 |
| 80 | 1.82.E-06 | 8.92.E-04 | 1.15.E-03 | 3.37.E-03 | 8.69.E-02 | 1.07.E+00 | 1.06.E+00 | 1.04.E+00 |
| 160 | 1.71.E-06 | 4.26.E-04 | 5.54.E-04 | 1.63.E-03 | 1.05.E-01 | 1.11.E+00 | 1.05.E+00 | 2.88.E-01 |
| 320 | 1.59.E-06 | 1.98.E-04 | 2.66.E-04 | 1.34.E-03 | | | | |

### 5.3 Application

We use $f(S) = \chi_{\{S=0\}}$ [39] as it is the simplest source term for penalizing sediment depletion $S = 0$. The problem is an HJBI equation with a discontinuous source. Alternatively, one may use a continuous regularized $f(S) = \kappa^{-1} \max\{S - \kappa, 0\}$ with a small $\kappa > 0$ if there is a regularity concern. Both provide the same numerical solution if $\kappa < \Delta S$.

We focus on the following four parameters affecting the optimal controls: observation cost $o$, uncertainty aversion coefficient $\psi$, and parameters $(W, \bar{L})$ controlling the waiting times. Unless otherwise specified, we choose $c_0 = 20$, $c_1 = 60$, $o = 20$, $\psi = 0.0001$, $\bar{L} = 10$, and $W = 2$ (day) such that the impacts of the model parameters become visible. The other parameter values are the same as those used in the test case. The resolution is $N_Q = N_S = 320$. In this case, the degree-of-freedom is 1,030,410. We then set $\varepsilon = 10^{-8}$.

**Figure 4** shows the computed $\Phi$ for $l = 1$ as a demonstrative example, showing that it is non-monotone and does not contain spurious oscillations. The computed $\Phi$ seems to be bounded as required in **Propositions 2-4**. The computed effective Hamiltonian $H$ is 0.590, which satisfies the bound stated in **Proposition 3**. Thus, the obtained computational results suggest that the proposed scheme reasonably handles the HJBI equation.

We compare the effective Hamiltonian $H$ for the different values of the parameters $o$, $\psi$, and $(W, \bar{L})$ with $W\bar{L} = 20$ (day) (**Table 8**). Using a larger (smaller) $\bar{L}$ results in a larger (smaller) number of degrees of freedom. The computational results imply a monotonic dependence of $H$ with respect to each parameter. The computational results follow the theoretical bound obtained in **Proposition 3**, suggesting consistency between the theoretical and numerical solutions.

We focus on the optimal controls $(L^*, \eta^*)$ in this study because they represent long-term state-dependent replenishment strategies. **Figures 5** and **6** show the dependence of the controls on the observation cost $o$. A larger observation cost motivates the decision-maker to reduce the observation opportunities as intuitively understood, which is quantified by solving the HJBI equation. Instead, replenishment should be carried out for each observation under wider conditions for a larger cost. Interestingly, for all cases, relatively sparse observations should be chosen during the low-flow period near the threshold $Q=\hat{Q}$ of the occurrence of a sediment transport. Replenishment should not be carried out at high discharge according to the proposed model, owing to the fact that the sediment will be immediately flushed out at such discharges. In addition, in real cases, it is usually difficult to replenish the sediment during a flood. These characteristics are common to the other cases examined below, suggesting their robustness.

**Figures 7** shows that a dense observation is required for moderately large flow discharges, whereas a sparser observation is needed for relatively low discharges. This is considered to be due to the clustering nature of the streamflow discharge model where low discharges close to the lower-bound $\underline{Q}$ persist. Although the dependence of $L^*$ on $\psi$ is non-monotone, the functional form of $L^*$ within the phase space is similar among the cases examined here. Therefore, the decision-maker should be rationally inattentive during a low-flow regime. **Figure 8** shows that the sediment replenishment should be carried out under the model uncertainty as well, but with a smaller amount and in smaller areas within the state space under a stronger aversion to uncertainty. The area where the sediment should be replenished shrinks as $\psi$ increases.

Finally, the impacts of Erlangization are analyzed in **Figures 9** and **10**. Increasing parameter $\bar{L}$ while maintaining the constraint $\bar{L}W = 20$ (day) leads to a more gradual transition of the optimal control $L^*$ within the state space, suggesting its convergence toward a spatially continuous transition. Increasing the level of Erlangization reduces the effective Hamiltonian $H$ because more flexible strategies become possible. Although a mathematical analysis of convergence using Erlangization $\bar{L} \to +\infty$ is beyond the scope of this study, the computational results support the existence of a nontrivial limit. The profile of $L^*$ does not change drastically among different cases with $\bar{L} \geq 10$. The optimal replenishment strategies are not critically different among the different cases, suggesting their robustness against the observation strategy. Using the proposed model, the decision-maker can design a replenishment strategy and determine how rationally inattentive he/she should be in both adaptive and problem-dependent ways.

**Table 8.** Effective Hamiltonian $H$ with respect to $o$, $\psi$, and $(W, \bar{L})$.

| $o$ | $H$ | $\psi$ | $H$ | $\bar{L}$ | $H$ |
|---|---|---|---|---|---|
| 10 | 0.545 | 0.00005 | 0.589 | 1 | 0.935 |
| 15 | 0.568 | 0.0001 | 0.590 | 5 | 0.654 |
| 20 | 0.590 | 0.0005 | 0.597 | 10 | 0.590 |
| 25 | 0.609 | 0.001 | 0.607 | 15 | 0.562 |
| 30 | 0.627 | 0.005 | 0.675 | 20 | 0.546 |
| 35 | 0.644 | 0.01 | 0.745 | 30 | 0.530 |

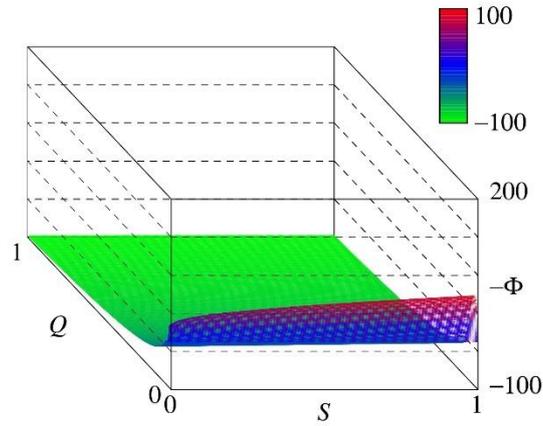

**Figure 4.** Computed potential $\Phi$ for $l = 1$ plotted as a 2-D surface within a 3-D space.

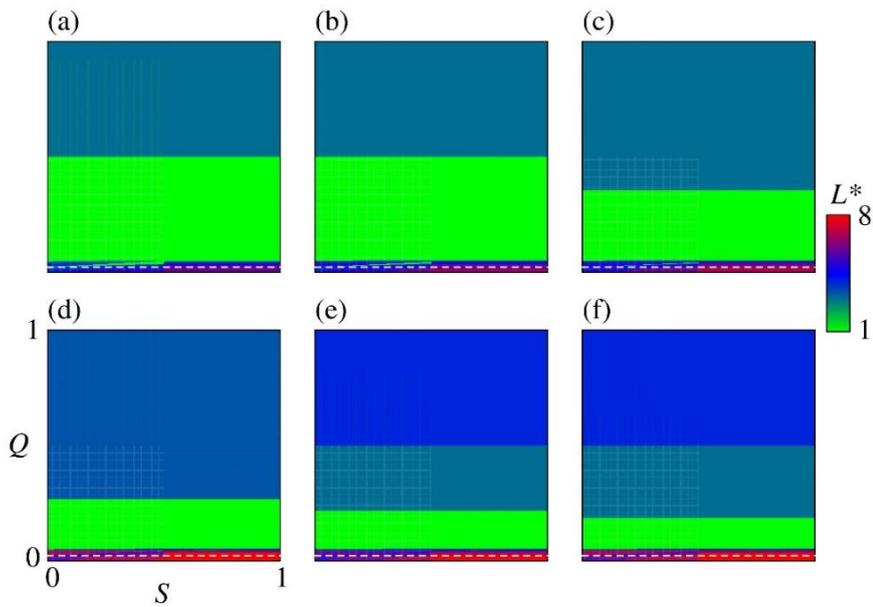

**Figure 5.** Optimal control $L^*$ for different $o$. The values of $o$ are (a) 10, (b) 15, (c) 20, (d) 25, (e) 30, and (f) 35. The dotted line represents the threshold $Q = \hat{Q}$.

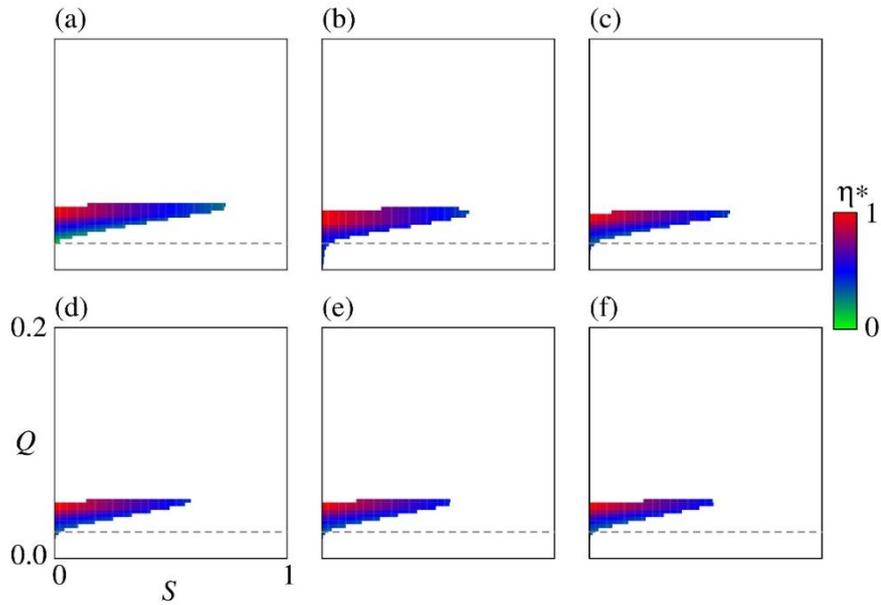

**Figure 6.** Optimal control $\eta^*$ for different $o$. The values of $o$ are (a) 10, (b) 15, (c) 20, (d) 25, (e) 30, and (f) 35. The dotted line represents the threshold $Q = \hat{Q}$, where $\eta^* = 0$ for $Q \geq 0.2$.

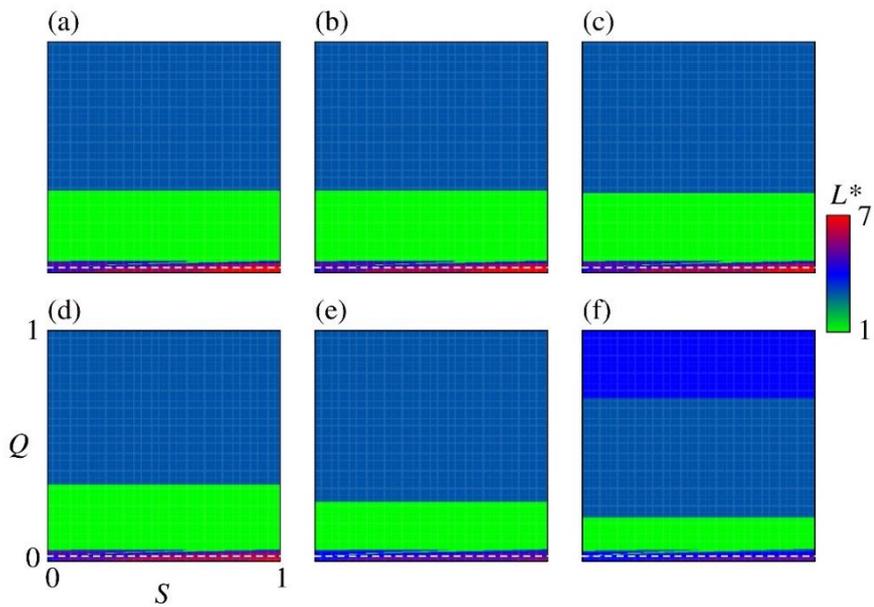

**Figure 7.** Optimal control $L^*$ for different $\psi$. The values of $\psi$ are (a) 0.00005, (b) 0.0001, (c) 0.0005, (d) 0.001, (e) 0.005, and (f) 0.01. The dotted line represents the threshold $Q = \hat{Q}$.

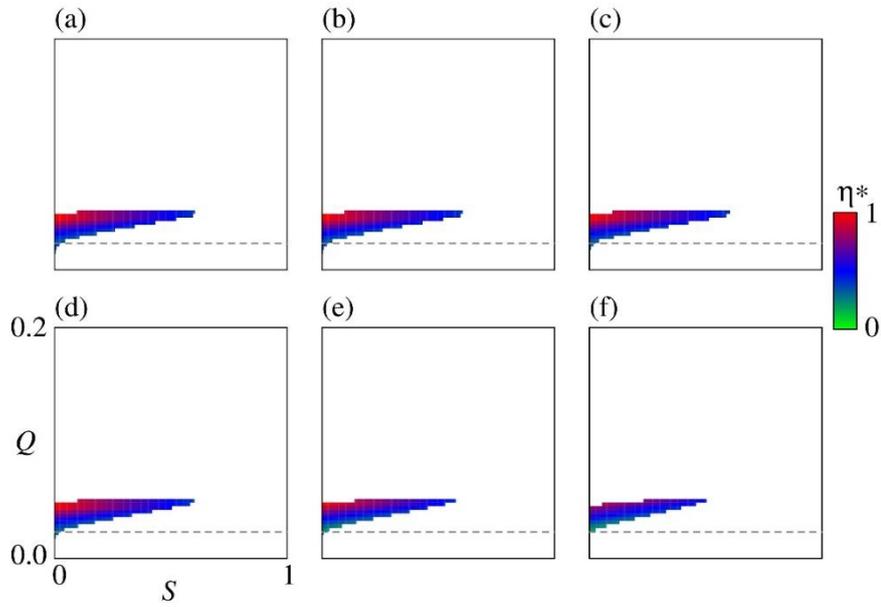

**Figure 8.** The optimal control $\eta^*$ for different values of $\psi$. The values of $\psi$ are (a) 0.00005, (b) 0.0001, (c) 0.0005, (d) 0.001, (e) 0.005, and (f) 0.01. The dotted line represents the threshold $Q = \hat{Q}$, where $\eta^* = 0$ for $Q \geq 0.2$.

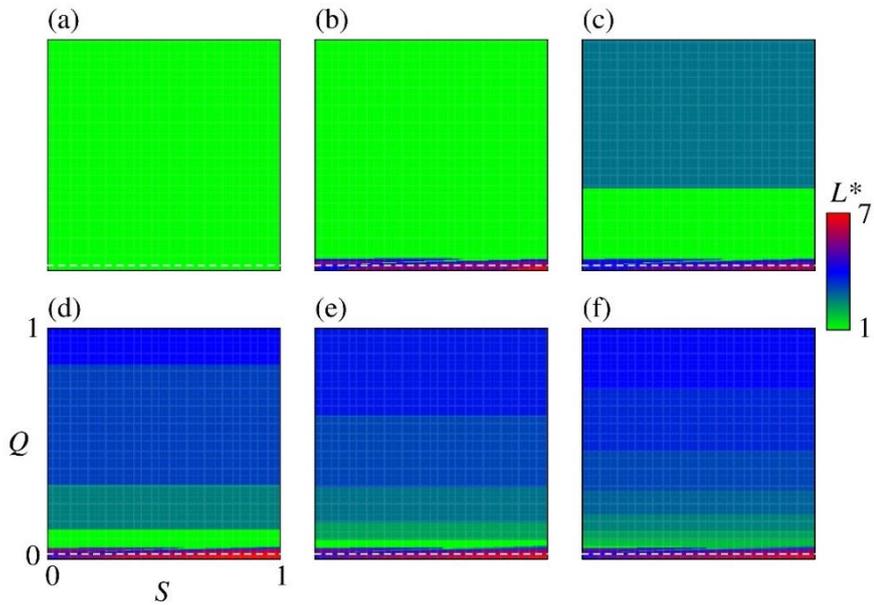

**Figure 9.** Optimal control $L^*$ for different $(W, \bar{L})$. The values of $\bar{L}$ are (a) 1, (b) 5, (c) 10, (d) 15, (e) 20, and (f) 30. The dotted line represents the threshold $Q = \hat{Q}$.

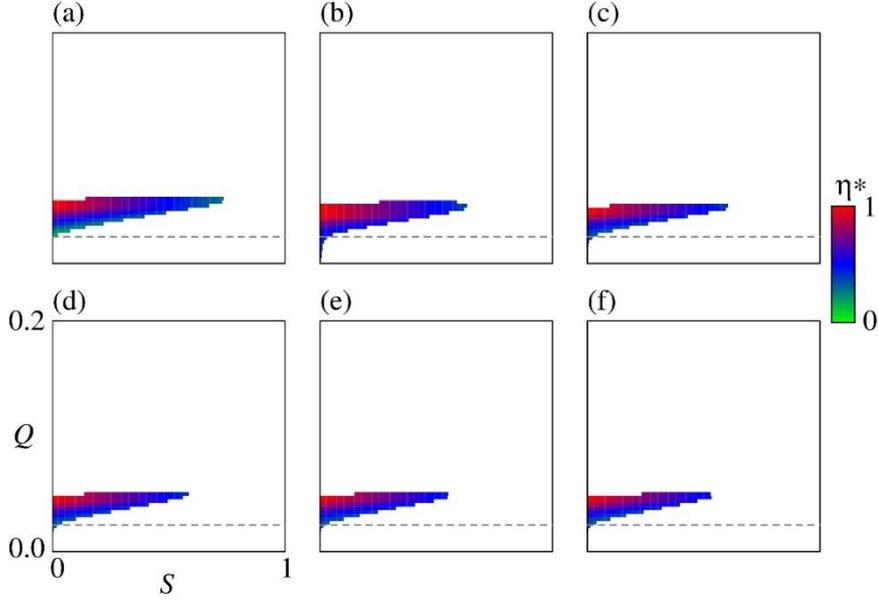

**Figure 10.** Optimal control $\eta^*$ for different $(W, \bar{L})$ with $W\bar{L} = 20$. The values of $\bar{L}$ are (a) 1, (b) 5, (c) 10, (d) 15, (e) 20, and (f) 30. The dotted line represents the threshold $Q = \hat{Q}$, where $\eta^* = 0$ for $Q \geq 0.2$.

# 6  Conclusions

We formulated and analyzed a stochastic control problem for sediment replenishment under uncertainty. Rational inattention as an adaptive observation/intervention strategy was effectively combined with dynamic programming under discrete observation/intervention. Erlangization plays a vital role in the efficient and adaptive mathematical modeling of the control problem. We discussed that finding the most cost-efficient replenishment strategy is reduced to solving an HJBI equation. Its closed-form solutions were not found, whereas the existence, uniqueness, and optimality results were derived in terms of the viscosity solution. Although we did not consider in this paper, the observation strategies can be restricted by replacing $\mathbb{M}\Phi$ at $l = 1$ using

$$\mathbb{M}_{\text{Lim}}\Phi(Q, S, 1) = \frac{1}{W}\Phi(Q, S, 1) - \frac{1}{W}\min_{i \in \Omega_L, u \in \Omega_S, i \geq \underline{i}}\left\{\Phi(Q, \min\{S + u, \bar{S}\}, i) + o + C(u)\right\}, \quad (50)$$

with a prescribed $\underline{i} \geq 1$ to consider sparser observations. The results of the optimality and uniqueness presented in this paper also apply to this case. Uncertainties in the sediment storage dynamics were also not considered but can be handled by extending the proposed framework. Detailed hydraulic experiments to quantify the range of parameters and their uncertainties will be conducted in the future.

The accuracy of the numerical scheme was verified based on a manufactured solution whose regularity can be tuned using a specific parameter. The numerical solutions of the scheme converged to the corresponding manufactured solutions with a convergence rate of no larger than 1. The demonstrative application of the model to a realistic case with identified parameter values suggested that an adaptive observation/intervention strategy based on rational inattention becomes optimal if the observation cost is considered. However, the regularity of the potentials solving HJBI and the related PIDEs arising in ergodic problems remains largely unknown. Our contribution therefore has theoretical restrictions. Nevertheless, to address this challenging issue, our contribution will become a case study under a unique problem setting. For example, one may reconsider the coupled streamflow-sediment-algae dynamics based on our framework [94]. An ergodic control problem can be considered as a small discount limit of a discounted problem or a large time limit of a finite-horizon problem. In such cases, dynamic programming principles, not directly utilized in this study, play central roles in a mathematical analysis. The gap between the verification and uniqueness results may be resolved using alternative approaches.

An alternative to the proposed approach is a continuous-state version of the imprecise transition probability model that does not require the benchmark and distorted models to be equivalent to each other [95]. This approach will provide more pessimistic results as the distorted models are explored in a wider space of probability measures. The existence and uniqueness of the optimal controls will become more difficult to analyze. We assumed no boundary singularity for an HJBI equation; however, for a problem with some boundary singularity, simply setting $\Phi = 0$ at a boundary point will not be possible and more sophisticated definitions of viscosity solutions will be necessary.

The use of a more generic stochastic process as the driving noise of eco-hydraulic dynamics, such as the genuinely nonlinear CBI processes [96], can be meaningful. This approach will become possible with high-frequency and long-term data, such as time series data with a frequency of several minutes and a length of several years. Finally, the Erlangization employed in this study can also be applied to other environmental management problems such as water quality management [97], which will also be an interesting research topic in the future. The construction of higher-order schemes is also an important topic for more accurate computations of HJBI equations. Applying upwind interpolation to the nonlocal term [98] is an implementable method. Mathematically, the convergence analysis of higher-order schemes is complicated because of the loss of monotonicity.

**Appendix A**

**A.1 Proof of Proposition 1**

We introduce a finite-horizon objective

$$J_T(Q,S;L,\eta,\phi) = \frac{1}{T}\mathbb{E}_{\mathbb{Q},0}\left[J_{1,T} + J_{2,T} + J_{3,T}\right] \tag{51}$$

such that

$$J(Q,S;L,\eta,\phi) = \limsup_{T\to+\infty} J_T(Q,S;L,\eta,\phi). \tag{52}$$

Here, $\mathbb{E}_{\mathbb{Q},0}[\cdot]$ is the expectation conditioned on $\mathcal{F}_0$ under $\mathbb{Q}$. This $J$ is a constant, assuming the existence of a unique invariant measure.

Fix $(Q_0, S_0, l_0) \in \Omega$ and choose some $(L, \eta, \phi) \in \mathcal{L} \times \mathcal{E} \times \mathcal{A}$. Based on the nature of Erlang, for each $L_{k+1} \in \mathbb{N}$ ($k = 0, 1, 2, ...$), we have a strictly increasing sequence $\{\varsigma_{k,m}\}_{m=0,2,3,...L_{k+1}}$ with $\tau_k = \varsigma_{k,0} < \varsigma_{k,1} < ... < \varsigma_{k,L_{k+1}-1} < \varsigma_{k,L_{k+1}} = \tau_{k+1}$ and i.i.d. $\varsigma_{k,m} - \varsigma_{k,m-1}$ ($1 \leq m \leq L_{k+1}$) following the exponential distribution with mean $W$. Define a càdlàg (right continuous with the left limit) process $(l_t)_{t \geq 0}$ by

$$l_t = \begin{cases} L_{k+1} & (t = \tau_k, \ k = 0, 1, 2, ...) \\ L_{k+1} - (m-1) & (\varsigma_{k,m-1} < t < \varsigma_{k,m}, \ k = 0, 1, 2, ..., \ m = 1, 2, 3, ..., L_{k+1}) \end{cases}. \quad (53)$$

Based on Itô's formula, for a generic regular $\Psi : \Omega \to \mathbb{R}$, we have

$$\begin{aligned}
\Psi(Q_T, S_T, l_T) - \Psi(Q_{\tau_k}, S_{\tau_k}, l_{\tau_k}) &= \Psi(Q_{\tau_k+}, S_{\tau_k+}, L_{k+1}) - \Psi(Q_{\tau_k}, S_{\tau_k}, 1) \\
&+ \sum_{\substack{l_t \neq l_{t-} \\ t \in (\tau_k, T)}} \left( \Psi(Q_t, S_t, l_t) - \Psi(Q_{t-}, S_{t-}, l_{t-}) \right) \\
&+ \int_{\tau_k}^T \left( -\rho(Q_t - \underline{Q}) \frac{\partial \Psi(Q_t, S_t, l_t)}{\partial Q} - F(Q_t, S_t) \frac{\partial \Psi(Q_t, S_t, l_t)}{\partial S} \right) \mathrm{d}t \quad (54) \\
&+ \int_{\tau_k}^T (Q_t + A) \int_0^\infty \left( \Psi(Q_t + z, S_t, l_t) - \Psi(Q_t, S_t, l_t) \right) v(\mathrm{d}z) \mathrm{d}t \\
&+ \int_{\tau_k}^T \int_0^\infty \int_0^{Q_{t-} + A} \left( \Psi(Q_{t-} + z, S_{t-}, l_{t-}) - \Psi(Q_{t-}, S_{t-}, l_{t-}) \right) \tilde{N}(\mathrm{d}z, \mathrm{d}u, \mathrm{d}t)
\end{aligned}$$

for $T \in (\tau_k, \tau_{k+1})$ ($k = 0, 1, 2, ...$). Each integration is applied for processes starting from $(Q_{\tau_k+}, S_{\tau_k+}, l_{\tau_k+})$. The second line represents the total $L_{k+1} - 1$ jumps of $\Psi(Q_t, S_t, l_t) - \Psi(Q_{t-}, S_{t-}, l_{t-})$ (see, (53)), which is dropped if $L_{k+1} = 1$ (exponential waiting time). We then have

$$\Psi(Q_{\tau_k+}, S_{\tau_k+}, L_{k+1}) = \Psi\left(Q_{\tau_k}, \min\{S_{\tau_k} + \eta_k, \bar{S}\}, L_{k+1}\right), \quad k = 1, 2, 3, .... \quad (55)$$

Fix $T > 0$, a small $\delta > 0$, and set $T_\delta = \min\{T, \inf\{\tau; Q_\tau > \delta^{-1}\}\}$. This $T_\delta$ is a.s. bounded and approaches $T$ as $\delta \to +0$. By definition, $Q_t$ is non-negative and bounded for $0 \leq t \leq T_\delta$, which is subsequently used in the following.

From (54) and (55), with $\Psi = \Phi$ we obtain

$$\Phi\left(Q_{T_\delta}, S_{T_\delta}, l_{T_\delta}\right) - \Phi(Q_0, S_0, l_0)$$

$$= \sum_{k=0}^{K} \left( \Phi\left(Q_{\tau_k+}, S_{\tau_k+}, l_{\tau_k+}\right) - \Phi\left(Q_{\tau_k}, S_{\tau_k}, l_{\tau_k}\right) \right)$$

$$+ \sum_{k=0}^{K} \sum_{\substack{l_t \neq l_{t-} \\ t \in (\tau_k, \min\{\tau_{k+1}, T_\delta\})}} \left( \Phi(Q_t, S_t, l_t) - \Phi(Q_{t-}, S_{t-}, l_{t-}) \right) \quad (56)$$

$$+ \sum_{k=0}^{K} \int_{\tau_k}^{\min\{\tau_{k+1}, T_\delta\}} \left( -\rho(Q_t - \underline{Q}) \frac{\partial \Phi(Q_t, S_t, l_t)}{\partial Q} - F(Q_t, S_t) \frac{\partial \Phi(Q_t, S_t, l_t)}{\partial S} \right) \mathrm{d}t$$

$$+ \sum_{k=0}^{K} \int_{\tau_k}^{\min\{\tau_{k+1}, T_\delta\}} (Q_t + A) \int_0^\infty \left( \Phi(Q_t + z, S_t, l_t) - \Phi(Q_t, S_t, l_t) \right) \phi_t(z) \nu(\mathrm{d}z) \mathrm{d}t$$

$$+ \sum_{k=0}^{K} \int_{\tau_k}^{\min\{\tau_{k+1}, T_\delta\}} \int_0^\infty \int_0^{Q_{t-}+A} \left( \Phi(Q_{t-} + z, S_{t-}, l_{t-}) - \Phi(Q_{t-}, S_{t-}, l_{t-}) \right) \tilde{N}(\mathrm{d}z, \mathrm{d}u, \mathrm{d}t)$$

using the notation $(Q_{0+}, S_{0+}, l_{0+}) = (Q_0, S_0, l_0)$, where $K$ is a non-negative integer with $T_\delta \in (\tau_K, \tau_{K+1})$.

Because $(L, \eta) \in \mathcal{L} \times \mathcal{E}$ and $\phi \in \mathcal{A}$ are arbitrary,

$$\mathbb{E}_{\mathbb{Q},0} \left[ \sum_{k=0}^{K} \left( \Phi\left(Q_{\tau_k}, \min\{S_{\tau_k} + \eta_k, \bar{S}\}, L_{k+1}\right) - \Phi\left(Q_{\tau_k}, S_{\tau_k}, 1\right) \right) \\ + \sum_{k=0}^{K} \sum_{\substack{l_t \neq l_{t-} \\ t \in \min\{\tau_{k+1}, T_\delta\}}} \left( \Phi(Q_t, S_t, l_t) - \Phi(Q_{t-}, S_{t-}, l_{t-}) \right) \right]$$

$$= \mathbb{E}_{\mathbb{Q},0} \left[ \sum_{k=0}^{K} \left( \begin{array}{l} \left( \Phi\left(Q_{\tau_k}, \min\{S_{\tau_k} + \eta_k, \bar{S}\}, L_{k+1}\right) - \Phi\left(Q_{\tau_k}, S_{\tau_k}, 1\right) + o + C(\eta_k) \right) \\ + \sum_{\substack{l_t \neq l_{t-} \\ t \in \min\{\tau_{k+1}, T_\delta\}}} \left( \Phi(Q_t, S_t, l_t) - \Phi(Q_{t-}, S_{t-}, l_{t-}) \right) \end{array} \right) \\ - \sum_{k=0}^{K} (o + C(\eta_k)) \right] \quad (57)$$

$$\geq \mathbb{E}_{\mathbb{Q},0} \left[ -\sum_{k=0}^{K} \int_{\tau_k}^{\min\{\tau_{k+1}, T_\delta\}} \mathbb{M} \Phi(Q_t, S_t, l_t) \mathrm{d}t \right] - \mathbb{E}_{\mathbb{Q}} \left[ \sum_{k=0}^{K} (o + C(\eta_k)) \right]$$

with equality for (20), and

$$-(Q_t + A) \sup_{\phi(\cdot) > 0} \left\{ -\int_0^\infty \left\{ \left( \Phi(Q_t, S_t, l_t) - \Phi(Q_t + z, S_t, l_t) \right) \phi(z) + \frac{1}{\psi} \left( \phi(z) \ln \phi(z) - \phi(z) + 1 \right) \right\} \nu(\mathrm{d}z) \right\}$$

$$-(Q_t + A) \int_0^\infty \left( \Phi(Q_t + z, S_t, l_t) - \Phi(Q_t, S_t, l_t) \right) \phi_t(z) \nu(\mathrm{d}z) \quad (58)$$

$$\leq (Q_t + A) \int_0^\infty \frac{1}{\psi} \left( \phi_t(z) \ln \phi_t(z) - \phi_t(z) + 1 \right) \nu(\mathrm{d}z)$$

with equality for (21).

First, we chose $\phi = \phi^* \in \mathcal{A}$ of (21). The probability measure induced by $\phi^*$ is denoted as $\mathbb{Q}^*$. By using (17) in a point-wise manner, we obtain the inequality

$$\mathbb{E}_{\mathbb{Q}^*,0}\left[\Phi\left(Q_{T_\delta},S_{T_\delta},l_{T_\delta}\right)-\Phi\left(Q_0,S_0,l_0\right)\right]$$

$$\geq -\mathbb{E}_{\mathbb{Q}^*,0}\left[\sum_{k=0}^{K}\left(o+C(\eta_k)\right)\right]$$

$$+\mathbb{E}_{\mathbb{Q}^*,0}\left[\sum_{k=0}^{K}\int_{\tau_k}^{\min\{\tau_{k+1},T_\delta\}}\binom{h-f(S_t)}{+(Q_t+A)\int_0^\infty \frac{1}{\psi}\left(\phi_t^*(z)\ln\phi_t^*(z)-\phi_t^*(z)+1\right)v(\mathrm{d}z)}\mathrm{d}t\right], \quad (59)$$

$$=\mathbb{E}_{\mathbb{Q}^*,0}\left[\min\{\tau_{K+1},T_\delta\}\right]h$$

$$-\mathbb{E}_{\mathbb{Q}^*,0}\left[\sum_{k=0}^{K}\left\{\begin{array}{l}\left(o+C(\eta_k)\right)+\int_{\tau_k}^{\min\{\tau_{k+1},T\}}f(S_t)\mathrm{d}t\\ -\int_{\tau_k}^{\min\{\tau_{k+1},T_\delta\}}(Q_t+A)\int_0^\infty \frac{1}{\psi}\left(\phi_t^*(z)\ln\phi_t^*(z)-\phi_t^*(z)+1\right)v(\mathrm{d}z)\mathrm{d}t\end{array}\right\}\right]$$

and thus

$$\frac{\mathbb{E}_{\mathbb{Q}^*,0}\left[\min\{\tau_{K+1},T_\delta\}\right]}{T_\delta}h$$

$$\leq \frac{1}{T_\delta}\mathbb{E}_{\mathbb{Q}^*,0}\left[\sum_{k=0}^{K}\left\{\begin{array}{l}\left(o+C(\eta_k)\right)+\int_{\tau_k}^{\min\{\tau_{k+1},T_\delta\}}f(S_t)\mathrm{d}t\\ -\int_{\tau_k}^{\min\{\tau_{k+1},T_\delta\}}(Q_t+A)\int_0^\infty \frac{1}{\psi}\left(\phi_t^*(z)\ln\phi_t^*(z)-\phi_t^*(z)+1\right)v(\mathrm{d}z)\mathrm{d}t\end{array}\right\}\right]. \quad (60)$$

$$+\frac{1}{T_\delta}\mathbb{E}_{\mathbb{Q}^*,0}\left[\Phi\left(Q_{T_\delta},S_{T_\delta},l_{T_\delta}\right)-\Phi\left(Q_0,S_0,l_0\right)\right]$$

By $\tau_{K+1}<\infty$, sending $\delta\to+\infty$ to obtain $T_\delta\to T$ and further

$$\lim_{T\to+\infty}\frac{\mathbb{E}_{\mathbb{Q}^*,0}\left[\min\{\tau_{K+1},T\}\right]}{T}=1, \quad (61)$$

leading to

$$h\leq \limsup_{T\to+\infty}J_T\left(Q,S;L,\eta,\phi^*\right) \quad (62)$$

and hence

$$h\leq \limsup_{T\to+\infty}J_T\left(Q,S;L^*,\eta^*,\phi^*\right)=H. \quad (63)$$

Second, choosing $\left(L^*,\eta^*\right)\in\mathcal{L}\times\mathcal{E}$ of (20), we then have

$$\mathbb{E}_{\mathbb{Q},0}\left[\Phi\left(Q_{T_\delta},S_{T_\delta},l_{T_\delta}\right)-\Phi(Q_0,S_0,l_0)\right]$$

$$\leq -\mathbb{E}_{\mathbb{Q},0}\left[\sum_{k=0}^{K}\left(o+C(\eta_k^*)\right)\right]$$

$$+\mathbb{E}_{\mathbb{Q},0}\left[\sum_{k=0}^{K}\int_{\tau_k}^{\min\{\tau_{k+1},T_\delta\}}\left(\begin{array}{c}h-f(S_t)\\+(Q_t+A)\int_0^\infty\frac{1}{\psi}(\phi_t(z)\ln\phi_t(z)-\phi_t(z)+1)v(\mathrm{d}z)\end{array}\right)\mathrm{d}t\right] \quad (64)$$

$$=\mathbb{E}_{\mathbb{Q},0}\left[\min\{\tau_{K+1},T_\delta\}\right]h$$

$$-\mathbb{E}_{\mathbb{Q},0}\left[\sum_{k=0}^{K}\left\{\begin{array}{c}\left(o+C(\eta_k^*)\right)+\int_{\tau_k}^{\min\{\tau_{k+1},T_\delta\}}f(S_t)\mathrm{d}t\\-\int_{\tau_k}^{\min\{\tau_{k+1},T_\delta\}}(Q_t+A)\int_0^\infty\frac{1}{\psi}(\phi_t(z)\ln\phi_t(z)-\phi_t(z)+1)v(\mathrm{d}z)\mathrm{d}t\end{array}\right\}\right]$$

and thus

$$\frac{\mathbb{E}_{\mathbb{Q},0}\left[\min\{\tau_{K+1},T_\delta\}\right]}{T_\delta}h$$

$$\geq \frac{1}{T_\delta}\mathbb{E}_{\mathbb{Q},0}\left[\sum_{k=0}^{K}\left\{\begin{array}{c}\left(o+C(\eta_k^*)\right)+\int_{\tau_k}^{\min\{\tau_{k+1},T_\delta\}}f(S_t)\mathrm{d}t\\-\int_{\tau_k}^{\min\{\tau_{k+1},T_\delta\}}(Q_t+A)\int_0^\infty\frac{1}{\psi}(\phi_t(z)\ln\phi_t(z)-\phi_t(z)+1)v(\mathrm{d}z)\mathrm{d}t\end{array}\right\}\right]. \quad (65)$$

$$+\frac{1}{T_\delta}\mathbb{E}_{\mathbb{Q},0}\left[\Phi\left(Q_{T_\delta},S_{T_\delta},l_{T_\delta}\right)-\Phi(Q_0,S_0,l_0)\right]$$

Then, as in the sub-solution case, we obtain

$$h\geq \limsup_{T\to+\infty}J_T\left(Q,S;L^*,\eta^*,\phi\right) \quad (66)$$

and hence

$$h\geq \limsup_{T\to+\infty}J_T\left(Q,S;L^*,\eta^*,\phi^*\right)=H. \quad (67)$$

Consequently, combining (63) and (67) yields the desired result.

## A.2 Proof of Proposition 2

We first consider the uncertainty-neutral case ($\psi\to+0$), and thus the uncertainty-averse case ($\psi>0$). The difference between the two cases is the nonlinearity of the integral term.

Set $\beta>0$. Based on the boundedness and continuity of sub- and super-solutions, the quantity $\Phi_1(Q,S,l)-\Phi_2(Q,S,l)-2\beta Q^2$ is maximized for $\Omega$ at $S=0$ or $S>0$. First, assume that it is maximized at $S>0$. The case with $S=0$ is handled in a similar way. Without a loss of generality, we assume $\Phi_1(Q,0,l)-\Phi_2(Q,0,l)<\Phi_1(Q,S,l)-\Phi_2(Q,S,l)$ for $S>0$. If, on the contrary, $\Phi_1(Q,0,l)-\Phi_2(Q,0,l)=\Phi_1(Q,S,l)-\Phi_2(Q,S,l)$ for $S>0$, we can resort to the case in which $\Phi_1-\Phi_2$ is maximized for $S=0$. Set $U_\beta=\max_{(Q,S,l)\in\Omega}\{\Phi_1(Q,S,l)-\Phi_2(Q,S,l)-2\beta Q^2\}$. We write the maximizer of $U_\beta$ as $(\hat{Q},\hat{S},\hat{l})$. Later, we send $\beta \to +0$.

Consider the following auxiliary function $\Xi_{\varepsilon,\beta}: \Omega \times \Omega \to \mathbb{R}$ with $\varepsilon > 0$:

$$\Xi_{\varepsilon,\beta}(Q_1,S_1,l_1,Q_2,S_2,l_2) = \Phi_1(Q_1,S_1,l_1) - \Phi_2(Q_2,S_2,l_2)$$
$$- \frac{1}{2\varepsilon}\left\{(Q_1-Q_2)^2 + (S_1-S_2)^2 + (l_1-l_2)^2\right\} - \beta(Q_1^2 + Q_2^2). \quad (68)$$

Set $U_{\varepsilon,\beta} = \sup_{\Omega \times \Omega} \Xi_{\varepsilon,\beta}$. Because $\Phi_1$ and $\Phi_2$ are bounded, the supremum is attained, i.e., $U_{\varepsilon,\beta} = \max_{\Omega \times \Omega} \Xi_{\varepsilon,\beta}$.

The maximizer of $\Xi_{\varepsilon,\beta}$ is denoted by $(\hat{Q}_1,\hat{S}_1,\hat{l}_1,\hat{Q}_2,\hat{S}_2,\hat{l}_2)$. By invoking a standard argument of the doubling of the variables [77, 99] along a subsequence (hereafter, we choose this subsequence), we obtain

$$\frac{1}{2\varepsilon}\left\{(\hat{S}_1-\hat{S}_2)^2 + (\hat{Q}_1-\hat{Q}_2)^2\right\} \to 0 \text{ and } \hat{l}_1 - \hat{l}_2 \to 0 \text{ as } \varepsilon \to +0. \quad (69)$$

Because $\hat{l}_1, \hat{l}_2 \in \mathbb{N}$, for a sufficiently small $\varepsilon > 0$, we have $\hat{l}_1 = \hat{l}_2$. Hence, a common limit $(\hat{Q},\hat{S},\hat{l}) \in \Omega$ exists for $(\hat{Q}_i,\hat{S}_i,\hat{l}_i)$ ($i=1,2$). We have $\Xi_{\varepsilon,\beta}(\hat{Q}_1,\hat{S}_1,\hat{l}_1,\hat{Q}_2,\hat{S}_2,\hat{l}_2) \geq \Xi_{\varepsilon,\beta}(\hat{Q},\hat{S},\hat{l},\hat{Q},\hat{S},\hat{l})$ and $U_\beta \geq U_{\varepsilon,\beta}$, and thus

$$\infty > U_\beta \geq U_{\varepsilon,\beta} \geq \Phi_1(\hat{Q},\hat{S},\hat{l}) - \Phi_2(\hat{Q},\hat{S},\hat{l}) - 2\beta\hat{Q}^2 > \frac{1}{2}U_\beta \text{ and } \hat{S}_1, \hat{S}_2 > 0 \quad (70)$$

for a small $\varepsilon > 0$, which is assumed hereafter. Based on the first inequality in (70), we obtain

$$\beta\hat{Q} \to 0 \text{ as } \beta \to 0. \quad (71)$$

The second inequality of (70) is inferred from the assumption that $\Phi_1(Q,S,l) - \Phi_2(Q,S,l) - 2\beta Q^2$ is maximized with $S > 0$.

Owing to **Definition 3.1** and (25) we obtain

$$H_i + \rho(\hat{Q}_i - \underline{Q})\frac{\partial \varphi_i(\hat{Q}_i,\hat{S}_i,\hat{l}_i)}{\partial Q} + F(\hat{Q}_i,\hat{S}_i)\frac{\partial \varphi_i(\hat{Q}_i,\hat{S}_i,\hat{l}_i)}{\partial S}$$
$$+ (\hat{Q}_i + A)\int_0^\delta \left(\varphi_i(\hat{Q}_i,\hat{S}_i,\hat{l}_i) - \varphi_i(\hat{Q}_i + z,\hat{S}_i,\hat{l}_i)\right)v(\mathrm{d}z) \quad \text{with } i=1 \text{ (resp., } i=2\text{)} \quad (72)$$
$$+ (\hat{Q}_i + A)\int_\delta^\infty \left(\Phi_i(\hat{Q}_i,\hat{S}_i,\hat{l}_i) - \Phi_i(\hat{Q}_i + z,\hat{S}_i,\hat{l}_i)\right)v(\mathrm{d}z)$$
$$+ \mathbb{M}\Phi_i(\hat{Q}_i,\hat{S}_i,\hat{l}_i) - f(\hat{S}_i) \leq 0 \quad (\text{resp., } \geq 0)$$

for the test functions $\varphi_1(Q,S,l) = \frac{1}{2\varepsilon}\left\{(S-\hat{S}_2)^2 + (Q-\hat{Q}_2)^2\right\} + \beta Q^2 + const$ and

$\varphi_2(Q,S,l) = -\frac{1}{2\varepsilon}\left\{(S-\hat{S}_1)^2 + (Q-\hat{Q}_1)^2\right\} - \beta Q^2 + const$ with a small $\delta > 0$. From (72), we obtain

$$H_1 - H_2 + I^{(1)} + I^{(2)} + I^{(3)} + I^{(4)} \leq 0, \quad (73)$$

where

$$I^{(1)} = -f\left(\hat{S}_1\right) + f\left(\hat{S}_2\right) + \rho\left(\hat{Q}_1 - \underline{Q}\right)\frac{\partial \varphi_1\left(\hat{Q}_1, \hat{S}_1, \hat{l}_1\right)}{\partial Q} - \rho\left(\hat{Q}_2 - \underline{Q}\right)\frac{\partial \varphi_2\left(\hat{Q}_2, \hat{S}_2, \hat{l}_2\right)}{\partial Q}$$
$$+ F\left(\hat{Q}_1, \hat{S}_1\right)\frac{\partial \varphi_1\left(\hat{Q}_1, \hat{S}_1, \hat{l}_1\right)}{\partial S} - F\left(\hat{Q}_2, \hat{S}_2\right)\frac{\partial \varphi_2\left(\hat{Q}_2, \hat{S}_2, \hat{l}_2\right)}{\partial S} \quad (74)$$

$$I^{(2)} = \mathbb{M}\Phi_1\left(\hat{Q}_1, \hat{S}_1, \hat{l}_1\right) - \mathbb{M}\Phi_2\left(\hat{Q}_2, \hat{S}_2, \hat{l}_2\right), \quad (75)$$

$$I^{(3)} = \left(\hat{Q}_1 + A\right)\int_0^\delta \left(\varphi_1\left(\hat{Q}_1, \hat{S}_1, \hat{l}_1\right) - \varphi_1\left(\hat{Q}_1 + z, \hat{S}_1, \hat{l}_1\right)\right)v(\mathrm{d}z)$$
$$-\left(\hat{Q}_2 + A\right)\int_0^\delta \left(\varphi_2\left(\hat{Q}_2, \hat{S}_2, \hat{l}_2\right) - \varphi_2\left(\hat{Q}_2 + z, \hat{S}_2, \hat{l}_2\right)\right)v(\mathrm{d}z), \quad (76)$$

$$I^{(4)} = \left(\hat{Q}_1 + A\right)\int_\delta^\infty \left(\Phi_1\left(\hat{Q}_1, \hat{S}_1, \hat{l}_1\right) - \Phi_1\left(\hat{Q}_1 + z, \hat{S}_1, \hat{l}_1\right)\right)v(\mathrm{d}z)$$
$$-\left(\hat{Q}_2 + A\right)\int_\delta^\infty \left(\Phi_2\left(\hat{Q}_2, \hat{S}_2, \hat{l}_2\right) - \Phi_2\left(\hat{Q}_2 + z, \hat{S}_2, \hat{l}_2\right)\right)v(\mathrm{d}z). \quad (77)$$

We analyze the convergence of each term below.

For $I^{(1)}$, we proceed as

$$I^{(1)} = -f\left(\hat{S}_1\right) + f\left(\hat{S}_2\right) + \rho\left(\hat{Q}_1 - \underline{Q}\right)\left\{\frac{\hat{Q}_1 - \hat{Q}_2}{\varepsilon} + 2\beta Q_1\right\} - \rho\left(\hat{Q}_2 - \underline{Q}\right)\left(\frac{\hat{Q}_1 - \hat{Q}_2}{\varepsilon} - 2\beta Q_2\right)$$
$$+ F\left(\hat{Q}_1, \hat{S}_1\right)\frac{1}{\varepsilon}\left(\hat{S}_1 - \hat{S}_2\right) - F\left(\hat{Q}_2, \hat{S}_2\right)\frac{1}{\varepsilon}\left(\hat{S}_1 - \hat{S}_2\right)$$
$$= -f\left(\hat{S}_1\right) + f\left(\hat{S}_2\right) + \frac{1}{\varepsilon}\rho\left(\hat{Q}_1 - \hat{Q}_2\right)^2 + 2\beta\rho\left(\hat{Q}_1 - \underline{Q}\right)Q_1 + 2\beta\rho\left(\hat{Q}_2 - \underline{Q}\right)Q_2$$
$$+ \frac{1}{\varepsilon}\left(F\left(\hat{Q}_1, \hat{S}_1\right) - F\left(\hat{Q}_2, \hat{S}_2\right)\right)\left(\hat{S}_1 - \hat{S}_2\right)$$
$$= -f\left(\hat{S}_1\right) + f\left(\hat{S}_2\right) \quad (78)$$
$$+ \frac{1}{\varepsilon}\rho\left(\hat{Q}_1 - \hat{Q}_2\right)^2 + 2\beta\rho\left(\hat{Q}_1 - \underline{Q}\right)Q_1 + 2\beta\rho\left(\hat{Q}_2 - \underline{Q}\right)Q_2$$
$$+ \frac{1}{\varepsilon}\left(F\left(\hat{Q}_1, \hat{S}_1\right) - F\left(\hat{Q}_2, \hat{S}_1\right)\right)\left(\hat{S}_1 - \hat{S}_2\right) + \frac{1}{\varepsilon}\left(F\left(\hat{Q}_2, \hat{S}_1\right) - F\left(\hat{Q}_2, \hat{S}_2\right)\right)\left(\hat{S}_1 - \hat{S}_2\right)$$
$$\geq -f\left(\hat{S}_1\right) + f\left(\hat{S}_2\right) + \frac{1}{\varepsilon}\rho\left(\hat{Q}_1 - \hat{Q}_2\right)^2$$
$$+ 2\beta\rho\left(\hat{Q}_1 - \underline{Q}\right)Q_1 + 2\beta\rho\left(\hat{Q}_2 - \underline{Q}\right)Q_2 + \frac{1}{\varepsilon}\left(F\left(\hat{Q}_1, \hat{S}_1\right) - F\left(\hat{Q}_2, \hat{S}_1\right)\right)\left(\hat{S}_1 - \hat{S}_2\right)$$

owing to the nondecreasing nature of $F$. Owing to the Lipschitz continuity of $F$, we obtain

$$I^{(1)} \geq 2 \times 2\beta\rho\left(\hat{Q} - \underline{Q}\right)\hat{Q} \quad \text{as} \quad \varepsilon \to +0. \quad (79)$$

For $I^{(2)}$, we have

$$WI^{(2)} = W\mathbb{M}\Phi_1\left(\hat{Q}_1, \hat{S}_1, \hat{l}_1\right) - W\mathbb{M}\Phi_2\left(\hat{Q}_2, \hat{S}_2, \hat{l}_2\right)$$

$$= \Phi_1\left(\hat{Q}_1, \hat{S}_1, \hat{l}_1\right) - \Phi_2\left(\hat{Q}_2, \hat{S}_2, \hat{l}_2\right)$$

$$- \begin{cases} \Phi_1\left(\hat{Q}_1, \hat{S}_1, \hat{l}_1 - 1\right) & \left(2 \leq \hat{l}_1 \leq L\right) \\ \min_{(i,\eta)\in\Omega_L\times\Omega_S}\left\{\Phi_1\left(\hat{Q}_1, \min\{\hat{S}_1+\eta, \bar{S}\}, i\right) + o + C(\eta)\right\} & \left(\hat{l}_1 = 1\right) \end{cases}$$

$$+ \begin{cases} \Phi_2\left(\hat{Q}_2, \hat{S}_2, \hat{l}_2 - 1\right) & \left(2 \leq \hat{l}_2 \leq L\right) \\ \min_{(i,\eta)\in\Omega_L\times\Omega_S}\left\{\Phi_2\left(\hat{Q}_2, \min\{\hat{S}_2+\eta, \bar{S}\}, i\right) + o + C(\eta)\right\} & \left(\hat{l}_2 = 1\right) \end{cases}. \quad (80)$$

Because $\hat{l}_1 = \hat{l}_2$, if $2 \leq \hat{l}_1 \leq L$, then

$$WI^{(2)} = W\mathbb{M}\Phi_1\left(\hat{Q}_1, \hat{S}_1, \hat{l}_1\right) - W\mathbb{M}\Phi_2\left(\hat{Q}_2, \hat{S}_2, \hat{l}_2\right)$$
$$= \Phi_1\left(\hat{Q}_1, \hat{S}_1, \hat{l}_1\right) - \Phi_2\left(\hat{Q}_2, \hat{S}_2, \hat{l}_2\right) - \Phi_1\left(\hat{Q}_1, \hat{S}_1, \hat{l}_1 - 1\right) + \Phi_2\left(\hat{Q}_2, \hat{S}_2, \hat{l}_2 - 1\right). \quad (81)$$

By the maximizing property of $\Xi_{\varepsilon,\beta}$, we have

$$\Phi_1\left(\hat{Q}_1, \hat{S}_1, \hat{l}_1\right) - \Phi_2\left(\hat{Q}_2, \hat{S}_2, \hat{l}_2\right) - \frac{1}{2\varepsilon}\left\{\left(\hat{S}_1 - \hat{S}_2\right)^2 + \left(\hat{Q}_1 - \hat{Q}_2\right)^2\right\} - \beta\left(\hat{Q}_1^2 + \hat{Q}_2^2\right)$$

$$= \Xi_{\varepsilon,\beta}\left(\hat{Q}_1, \hat{S}_1, \hat{l}_1, \hat{Q}_2, \hat{S}_2, \hat{l}_2\right)$$

$$\geq \Xi_{\varepsilon,\beta}\left(\hat{Q}_1, \hat{S}_1, \hat{l}_1 - 1, \hat{Q}_2, \hat{S}_2, \hat{l}_2 - 1\right) \quad (82)$$

$$= \Phi_1\left(\hat{Q}_1, \hat{S}_1, \hat{l}_1 - 1\right) - \Phi_2\left(\hat{Q}_2, \hat{S}_2, \hat{l}_2 - 1\right) - \frac{1}{2\varepsilon}\left\{\left(\hat{S}_1 - \hat{S}_2\right)^2 + \left(\hat{Q}_1 - \hat{Q}_2\right)^2\right\} - \beta\left(\hat{Q}_1^2 + \hat{Q}_2^2\right)$$

and thus

$$\Phi_1\left(\hat{Q}_1, \hat{S}_1, \hat{l}_1\right) - \Phi_2\left(\hat{Q}_2, \hat{S}_2, \hat{l}_2\right) - \Phi_1\left(\hat{Q}_1, \hat{S}_1, \hat{l}_1 - 1\right) + \Phi_2\left(\hat{Q}_2, \hat{S}_2, \hat{l}_2 - 1\right) \geq 0. \quad (83)$$

By (81) and (83), we obtain $I^{(2)} \geq 0$ for $\varepsilon, \beta > 0$.

A minimizer of $\mathbb{M}\Phi_i$ is denoted as $\left(l^{(i)}, \eta^{(i)}\right) \in \Omega_L \times \Omega_S$ ($i = 1, 2$). If $\hat{l}_1\left(=\hat{l}_2\right) = 1$, then

$$WI^{(2)} = \Phi_1\left(\hat{Q}_1, \hat{S}_1, 1\right) - \Phi_2\left(\hat{Q}_2, \hat{S}_2, 1\right)$$

$$- \min_{(l,\eta)\in\Omega_L\times\Omega_S}\left\{\Phi_1\left(\hat{Q}_1, \min\{\hat{S}_1+\eta, \bar{S}\}, i\right) + o + C(\eta)\right\}$$

$$+ \min_{(l,\eta)\in\Omega_L\times\Omega_S}\left\{\Phi_2\left(\hat{Q}_2, \min\{\hat{S}_2+\eta, \bar{S}\}, i\right) + o + C(\eta)\right\}$$

$$= \Phi_1\left(\hat{Q}_1, \hat{S}_1, 1\right) - \Phi_2\left(\hat{Q}_2, \hat{S}_2, 1\right) \quad (84)$$

$$+ \left\{\Phi_2\left(\hat{Q}_2, \min\{\hat{S}_2+\eta^{(2)}, \bar{S}\}, l^{(2)}\right) + C\left(\eta^{(2)}\right)\right\} - \left\{\Phi_1\left(\hat{Q}_1, \min\{\hat{S}_1+\eta^{(1)}, \bar{S}\}, l^{(1)}\right) + C\left(\eta^{(1)}\right)\right\}$$

$$\leq \Phi_1\left(\hat{Q}_1, \hat{S}_1, 1\right) - \Phi_2\left(\hat{Q}_2, \hat{S}_2, 1\right)$$

$$+ \left\{\Phi_2\left(\hat{Q}_2, \min\{\hat{S}_2+\eta^{(1)}, \bar{S}\}, l^{(1)}\right) + C\left(\eta^{(1)}\right)\right\} - \left\{\Phi_1\left(\hat{Q}_1, \min\{\hat{S}_1+\eta^{(1)}, \bar{S}\}, l^{(1)}\right) + C\left(\eta^{(1)}\right)\right\}$$

$$= \Phi_1\left(\hat{Q}_1, \hat{S}_1, 1\right) - \Phi_2\left(\hat{Q}_2, \hat{S}_2, 1\right) + \Phi_2\left(\hat{Q}_2, \min\{\hat{S}_2+\eta^{(1)}, \bar{S}\}, l^{(1)}\right) - \Phi_1\left(\hat{Q}_1, \min\{\hat{S}_1+\eta^{(1)}, \bar{S}\}, l^{(1)}\right)$$

Again, by the maximizing property of $\Xi_{\varepsilon,\beta}$, we have

$$\Phi_1(\hat{Q}_1,\hat{S}_1,\hat{I}_1) - \Phi_2(\hat{Q}_2,\hat{S}_2,\hat{I}_2) - \frac{1}{2\varepsilon}\left\{(\hat{S}_1-\hat{S}_2)^2+(\hat{Q}_1-\hat{Q}_2)^2\right\} - \beta(\hat{Q}_1^2+\hat{Q}_2^2)$$

$$= \Xi_{\varepsilon,\beta}(\hat{Q}_1,\hat{S}_1,\hat{I}_1,\hat{Q}_2,\hat{S}_2,\hat{I}_2)$$

$$\geq \Xi_{\varepsilon,\beta}\left(\hat{Q}_1, \min\{\hat{S}_1+\eta^{(1)},\overline{S}\}, I^{(1)}, \hat{Q}_2, \min\{\hat{S}_2+\eta^{(1)},\overline{S}\}, I^{(1)}\right) \qquad (85)$$

$$= \Phi_1\left(\hat{Q}_1, \min\{\hat{S}_1+\eta^{(1)},\overline{S}\}, I^{(1)}\right) - \Phi_2\left(\hat{Q}_2, \min\{\hat{S}_2+\eta^{(1)},\overline{S}\}, I^{(1)}\right)$$

$$-\frac{1}{2\varepsilon}\left\{\left(\min\{\hat{S}_1+\eta^{(1)},\overline{S}\}-\min\{\hat{S}_2+\eta^{(1)},\overline{S}\}\right)^2+(\hat{Q}_1-\hat{Q}_2)^2\right\} - \beta(\hat{Q}_1^2+\hat{Q}_2^2)$$

and thus

$$\Phi_1(\hat{Q}_1,\hat{S}_1,\hat{I}_1) - \Phi_2(\hat{Q}_2,\hat{S}_2,\hat{I}_2) + \Phi_2\left(\hat{Q}_2, \min\{\hat{S}_2+\eta^{(1)},\overline{S}\}, I^{(1)}\right) - \Phi_1\left(\hat{Q}_1, \min\{\hat{S}_1+\eta^{(1)},\overline{S}\}, I^{(1)}\right)$$

$$\geq -\frac{1}{2\varepsilon}\left(\min\{\hat{S}_1+\eta^{(1)},\overline{S}\}-\min\{\hat{S}_2+\eta^{(1)},\overline{S}\}\right)^2 + \frac{1}{2\varepsilon}(\hat{S}_1-\hat{S}_2)^2 \qquad (86)$$

An elementary calculation shows

$$\left|\min\{\hat{S}_1+\eta^{(1)},\overline{S}\}-\min\{\hat{S}_2+\eta^{(1)},\overline{S}\}\right| \leq \left|\hat{S}_1+\eta^{(1)}-(\hat{S}_2+\eta^{(1)})\right| = \left|\hat{S}_1-\hat{S}_2\right|. \qquad (87)$$

Then, (86) gives

$$\Phi_1(\hat{Q}_1,\hat{S}_1,\hat{I}_1) - \Phi_2(\hat{Q}_2,\hat{S}_2,\hat{I}_2) + \Phi_2\left(\hat{Q}_2, \min\{\hat{S}_2+\eta^{(1)},\overline{S}\}, I^{(1)}\right) - \Phi_1\left(\hat{Q}_1, \min\{\hat{S}_1+\eta^{(1)},\overline{S}\}, I^{(1)}\right)$$

$$\geq -\frac{1}{2\varepsilon}(\hat{S}_1-\hat{S}_2)^2 + \frac{1}{2\varepsilon}(\hat{S}_1-\hat{S}_2)^2 \qquad (88)$$

$$\geq 0$$

We again obtain $I^{(2)} \geq 0$ for $\varepsilon, \beta > 0$.

For $I^{(3)}$, we proceed as

$$I^{(3)} = (\hat{Q}_1+A)\int_0^\delta \left(\frac{1}{2\varepsilon}(\hat{Q}_1-\hat{Q}_2)^2+\beta\hat{Q}_1^2 - \frac{1}{2\varepsilon}(\hat{Q}_1+z-\hat{Q}_2)^2 - \beta(\hat{Q}_1+z)^2\right)v(\mathrm{d}z)$$

$$-(\hat{Q}_2+A)\int_0^\delta\left(-\frac{1}{2\varepsilon}(\hat{Q}_2-\hat{Q}_1)^2-\beta\hat{Q}_2^2+\frac{1}{2\varepsilon}(\hat{Q}_2+z-\hat{Q}_1)^2+\beta(\hat{Q}_2+z)^2\right)v(\mathrm{d}z)$$

$$= (\hat{Q}_1+A)\int_0^\delta\left(\frac{1}{2\varepsilon}(2\hat{Q}_1-2\hat{Q}_2+z)(-z)+\beta\hat{Q}_1^2-\beta(\hat{Q}_1+z)^2\right)v(\mathrm{d}z)$$

$$-(\hat{Q}_2+A)\int_0^\delta\left(\frac{1}{2\varepsilon}(2\hat{Q}_2-2\hat{Q}_1+z)z+\beta(\hat{Q}_2+z)^2-\beta\hat{Q}_2^2\right)v(\mathrm{d}z) \qquad , (89)$$

$$= -\beta(\hat{Q}_1+A)\int_0^\delta\left((\hat{Q}_1+z)^2-\hat{Q}_1^2\right)v(\mathrm{d}z) - \beta(\hat{Q}_2+A)\int_0^\delta\left((\hat{Q}_2+z)^2-\hat{Q}_2^2\right)v(\mathrm{d}z)$$

$$+(\hat{Q}_1+A)\int_0^\delta\frac{1}{2\varepsilon}(2\hat{Q}_2-2\hat{Q}_1-z)zv(\mathrm{d}z)-(\hat{Q}_2+A)\int_0^\delta\frac{1}{2\varepsilon}(2\hat{Q}_2-2\hat{Q}_1+z)zv(\mathrm{d}z)$$

$$= -\sum_{m=1}^2\sum_{j=0}^1(\hat{Q}_m+A)\binom{2}{j}\beta\hat{Q}_m^j\left(\int_0^\delta z^{2-j}v(\mathrm{d}z)\right) - \frac{1}{\varepsilon}(\hat{Q}_2-\hat{Q}_1)^2\int_0^\delta zv(\mathrm{d}z) - \left(\frac{\hat{Q}_1+\hat{Q}_2}{2}+A\right)\frac{1}{\varepsilon}\int_0^\delta z^2v(\mathrm{d}z)$$

Namely,

$$I^{(3)} = -\left(\frac{\hat{Q}_1+\hat{Q}_2}{2}+A\right)\frac{1}{\varepsilon}\int_0^\delta z^2v(\mathrm{d}z) - \frac{(\hat{Q}_2-\hat{Q}_1)^2}{\varepsilon}\int_0^\delta zv(\mathrm{d}z) - \sum_{m=1}^2\sum_{j=0}^1(\hat{Q}_m+A)\binom{2}{j}\beta\hat{Q}_m^j\left(\int_0^\delta z^{2-j}v(\mathrm{d}z)\right). \qquad (90)$$

Choose parameterization $\delta = \delta(\varepsilon) = \varepsilon^{\gamma}$ with $\gamma > (1-\alpha)^{-1}$. Then, on the last term of (90), we obtain

$$\frac{1}{\varepsilon}\int_0^{\delta} z^{2-j} v(\mathrm{d}z) \leq \frac{a}{\varepsilon}\int_0^{\delta}\frac{z^{2-j}}{z^{1+\alpha}}\mathrm{d}z = \frac{a}{2-j-\alpha}\frac{\delta^{2-j-\alpha}}{\varepsilon} = \frac{a}{2-j-\alpha}\varepsilon^{(2-j-\alpha)\gamma - 1} \to 0 \text{ as } \varepsilon \to +0 \quad (91)$$

because $j = 0$ or $1$. By (71), we obtain

$$I^{(3)} \to 0 \text{ as } \varepsilon \to +0 \text{ for each } \beta > 0. \quad (92)$$

For $I^{(4)}$, for $z \geq 0$ we have

$$\Phi_1(\hat{Q}_1,\hat{S}_1,\hat{I}_1) - \Phi_2(\hat{Q}_2,\hat{S}_2,\hat{I}_2) - \frac{1}{2\varepsilon}\left\{(\hat{S}_1 - \hat{S}_2)^2 + (\hat{Q}_1 - \hat{Q}_2)^2\right\} - \beta(\hat{Q}_1^2 + \hat{Q}_2^2)$$
$$= \Xi_{\varepsilon,\beta}(\hat{Q}_1,\hat{S}_1,\hat{I}_1,\hat{Q}_2,\hat{S}_2,\hat{I}_2)$$
$$\geq \Xi_{\varepsilon,\beta}(\hat{Q}_1 + z,\hat{S}_1,\hat{I}_1,\hat{Q}_2 + z,\hat{S}_2,\hat{I}_2) \quad , \quad (93)$$
$$= \Phi_1(\hat{Q}_1 + z,\hat{S}_1,\hat{I}_1) - \Phi_2(\hat{Q}_2 + z,\hat{S}_2,\hat{I}_2) - \frac{1}{2\varepsilon}\left\{(\hat{S}_1 - \hat{S}_2)^2 + (\hat{Q}_1 + z - (\hat{Q}_2 + z))^2\right\}$$
$$- \beta\left\{(\hat{Q}_1 + z)^2 + (\hat{Q}_2 + z)^2\right\}$$

and thus

$$\Phi_1(\hat{Q}_1,\hat{S}_1,\hat{I}_1) - \Phi_2(\hat{Q}_2,\hat{S}_2,\hat{I}_2) - \Phi_1(\hat{Q}_1 + z,\hat{S}_1,\hat{I}_1) + \Phi_2(\hat{Q}_2 + z,\hat{S}_2,\hat{I}_2)$$
$$\geq \beta(\hat{Q}_1^2 - (\hat{Q}_1 + z)^2) + \beta(\hat{Q}_2^2 - (\hat{Q}_2 + z)^2) \quad . \quad (94)$$

Hence, we can proceed as follows:

$$I^{(4)} = (\hat{Q}_1 + A)\int_{\delta}^{\infty}(\Phi_1(\hat{Q}_1,\hat{S}_1,\hat{I}_1) - \Phi_1(\hat{Q}_1 + z,\hat{S}_1,\hat{I}_1))v(\mathrm{d}z)$$
$$- (\hat{Q}_2 + A)\int_{\delta}^{\infty}(\Phi_2(\hat{Q}_2,\hat{S}_2,\hat{I}_2) - \Phi_2(\hat{Q}_2 + z,\hat{S}_2,\hat{I}_2))v(\mathrm{d}z)$$
$$= (\hat{Q}_1 - \hat{Q}_2)\int_{\delta}^{\infty}(\Phi_1(\hat{Q}_1,\hat{S}_1,\hat{I}_1) - \Phi_1(\hat{Q}_1 + z,\hat{S}_1,\hat{I}_1))v(\mathrm{d}z)$$
$$+ (\hat{Q}_2 + A)\int_{\delta}^{\infty}(\Phi_1(\hat{Q}_1,\hat{S}_1,\hat{I}_1) - \Phi_1(\hat{Q}_1 + z,\hat{S}_1,\hat{I}_1))v(\mathrm{d}z) \quad . \quad (95)$$
$$- (\hat{Q}_2 + A)\int_{\delta}^{\infty}(\Phi_2(\hat{Q}_2,\hat{S}_2,\hat{I}_2) - \Phi_2(\hat{Q}_2 + z,\hat{S}_2,\hat{I}_2))v(\mathrm{d}z)$$
$$\geq (\hat{Q}_1 - \hat{Q}_2)\int_{\delta}^{\infty}(\Phi_1(\hat{Q}_1,\hat{S}_1,\hat{I}_1) - \Phi_1(\hat{Q}_1 + z,\hat{S}_1,\hat{I}_1))v(\mathrm{d}z)$$
$$- (\hat{Q}_2 + A)\int_{\delta}^{\infty}\beta\left((\hat{Q}_1 + z)^2 - \hat{Q}_1^2 + (\hat{Q}_2 + z)^2 - \hat{Q}_2^2\right)v(\mathrm{d}z)$$

By the Hölder continuity of $\Phi_1$ with $\theta > \alpha$ and (69), the first term is evaluated as

$$\left|(\hat{Q}_1 - \hat{Q}_2)\int_{\delta}^{\infty}(\Phi_1(\hat{Q}_1,\hat{S}_1,\hat{I}_1) - \Phi_1(\hat{Q}_1 + z,\hat{S}_1,\hat{I}_1))v(\mathrm{d}z)\right| \leq |\hat{Q}_1 - \hat{Q}_2|C^{(0)}\int_{\delta}^{\infty}z^{\theta}v(\mathrm{d}z)$$
$$\leq |\hat{Q}_1 - \hat{Q}_2|C^{(0)}\int_0^{\infty}z^{\theta}v(\mathrm{d}z) \to 0 \quad , \quad \varepsilon \to +0. \quad (96)$$

On the last integral of (95), we have

$$-\left(\hat{Q}_2+A\right)\int_\delta^\infty \beta\left(\left(\hat{Q}_1+z\right)^2-\hat{Q}_1^2+\left(\hat{Q}_2+z\right)^2-\hat{Q}_2^2\right)v(\mathrm{d}z)$$

$$\to -2\beta\left(\hat{Q}+A\right)\int_0^\infty\left(\left(\hat{Q}+z\right)^2-\hat{Q}^2\right)v(\mathrm{d}z)$$

$$=-2\beta\left(\hat{Q}+A\right)\left(2\hat{Q}\int_0^\infty zv(\mathrm{d}z)+\int_0^\infty z^2 v(\mathrm{d}z)\right)$$

$$=-\beta\left(\hat{Q}+A\right)\left(4\hat{Q}M_1+2M_2\right)$$

$(m=1,2)$ as $\varepsilon\to +0$. (97)

Now, we use (79) and (97). By $\rho-M_1>0$, we obtain

$$\begin{aligned}&\beta\left\{4\rho\hat{Q}^2-4\rho\underline{Q}\hat{Q}-\left(\hat{Q}+A\right)\left(4\hat{Q}M_1+2M_2\right)\right\}\\&=2\beta\left\{2\rho\hat{Q}^2-2\rho\underline{Q}\hat{Q}-\left(\hat{Q}+A\right)\left(2\hat{Q}M_1+M_2\right)\right\}\\&\geq 2\beta\left\{2(\rho-M_1)\hat{Q}^2-2\rho\underline{Q}\hat{Q}-M_2\hat{Q}-A\left(2\hat{Q}M_1+M_2\right)\right\},\\&\geq \beta\left\{-2\rho\underline{Q}\hat{Q}-M_2\hat{Q}-A\left(2\hat{Q}M_1+M_2\right)\right\}\\&\to 0\text{ as }\beta\to +0\end{aligned}$$ (98)

where we use (71) in the last line. Combining (79) and (95)–(98) yields

$$I^{(1)}+I^{(4)}\geq 0\text{ as }\varepsilon\to +0\text{ and then }\beta\to +0.$$ (99)

Now, we consider the limit of (73). The calculation above shows that each $I^{(i)}\geq 0$ by letting $\varepsilon\to +0$ and then $\beta\to +0$ with $\delta=\varepsilon^\gamma$ ($\gamma>(1-\alpha)^{-1}$). Under this limit, we obtain $I^{(1)}+I^{(2)}+I^{(3)}+I^{(4)}\geq 0$, and hence $H_1-H_2\leq 0$ by (73), which is the desired result.

For an uncertainty-averse case ($\psi>0$), (76) and (77) are replaced with

$$I^{(3)}=\frac{\hat{Q}_1+A}{\psi}\int_0^\delta\left\{1-e^{-\psi\left(\varphi_1(\hat{Q}_1,\hat{S}_1,\hat{I}_1)-\varphi_1(\hat{Q}_1+z,\hat{S}_1,\hat{I}_1)\right)}\right\}v(\mathrm{d}z)-\frac{\hat{Q}_2+A}{\psi}\int_0^\delta\left\{1-e^{-\psi\left(\varphi_2(\hat{Q}_2,\hat{S}_2,\hat{I}_2)-\varphi_2(\hat{Q}_2+z,\hat{S}_2,\hat{I}_2)\right)}\right\}v(\mathrm{d}z),$$ (100)

$$I^{(4)}=\frac{\hat{Q}_1+A}{\psi}\int_\delta^\infty\left\{1-e^{-\psi\left(\Phi_1(\hat{Q}_1,\hat{S}_1,\hat{I}_1)-\Phi_1(\hat{Q}_1+z,\hat{S}_1,\hat{I}_1)\right)}\right\}v(\mathrm{d}z)-\frac{\hat{Q}_2+A}{\psi}\int_\delta^\infty\left\{1-e^{-\psi\left(\Phi_2(\hat{Q}_2,\hat{S}_2,\hat{I}_2)-\Phi_2(\hat{Q}_2+z,\hat{S}_2,\hat{I}_2)\right)}\right\}v(\mathrm{d}z).$$ (101)

The convergence of this $I^{(3)}$ as $\varepsilon\to +0$ is concluded as follows: First, we have

$$\begin{aligned}I^{(3)}=&\frac{\hat{Q}_1+A}{\psi}\int_0^\delta\left\{1-e^{\psi\left(\frac{1}{2\varepsilon}(2\hat{Q}_1-2\hat{Q}_2+z)z+\beta\hat{Q}_1^2-\beta(\hat{Q}_1+z)^2\right)}\right\}v(\mathrm{d}z)\\&-\frac{\hat{Q}_2+A}{\psi}\int_0^\delta\left\{1-e^{-\psi\left(\frac{1}{2\varepsilon}(2\hat{Q}_2-2\hat{Q}_1+z)z-\beta\hat{Q}_2^2+\beta(\hat{Q}_2+z)^2\right)}\right\}v(\mathrm{d}z)\end{aligned}.$$ (102)

By choosing $\delta=\delta(\varepsilon)=\varepsilon^\gamma$ with $\gamma>\max\left\{2^{-1},(1-\alpha)^{-1}\right\}$, for $0\leq z\leq\delta$ we have

$$\left|\frac{1}{2\varepsilon}\left(2\hat{Q}_1-2\hat{Q}_2+z\right)z\right|,\left|\left(\frac{1}{2\varepsilon}\left(2\hat{Q}_2-2\hat{Q}_1+z\right)z\right)\right|\leq\frac{\left|\hat{Q}_1-\hat{Q}_2\right|}{\varepsilon}z+\frac{z^2}{2\varepsilon}\leq\frac{\left|\hat{Q}_1-\hat{Q}_2\right|}{\varepsilon^{1/2}}\varepsilon^{\gamma-1/2}+\frac{\varepsilon^{2\gamma-1}}{2}\to 0$$ (103)

as $\varepsilon\to +0$. Furthermore, for $0\leq z\leq\delta(\varepsilon)$, we have $\left(\hat{Q}_m+z\right)^2-\hat{Q}_m^2\to 0$ $\varepsilon\to +0$. The integrands of

(102) are uniformly bound for $\delta = \delta(\varepsilon) = \varepsilon^{\gamma}$ and a small $\varepsilon$, and hence $\lim_{\varepsilon \to +0} I^{(3)} \geq 0$.

On $I^{(4)}$, we first have

$$I^{(4)} = \frac{\hat{Q}_1 + A}{\psi} \int_\delta^\infty \left\{1 - e^{-\psi\left(\Phi_1(\hat{Q}_1, \hat{S}_1, \hat{I}_1) - \Phi_1(\hat{Q}_1 + z, \hat{S}_1, \hat{I}_1)\right)}\right\} v(\mathrm{d}z) - \frac{\hat{Q}_2 + A}{\psi} \int_\delta^\infty \left\{1 - e^{-\psi\left(\Phi_2(\hat{Q}_2, \hat{S}_2, \hat{I}_2) - \Phi_2(\hat{Q}_2 + z, \hat{S}_2, \hat{I}_2)\right)}\right\} v(\mathrm{d}z)$$

$$= \frac{\hat{Q}_1 - \hat{Q}_2}{\psi} \int_\delta^\infty \left\{1 - e^{-\psi\left(\Phi_1(\hat{Q}_1, \hat{S}_1, \hat{I}_1) - \Phi_1(\hat{Q}_1 + z, \hat{S}_1, \hat{I}_1)\right)}\right\} v(\mathrm{d}z) \qquad (104)$$

$$+ \frac{\hat{Q}_2 + A}{\psi} \int_\delta^\infty \left\{e^{-\psi\left(\Phi_2(\hat{Q}_2, \hat{S}_2, \hat{I}_2) - \Phi_2(\hat{Q}_2 + z, \hat{S}_2, \hat{I}_2)\right)} - e^{-\psi\left(\Phi_1(\hat{Q}_1, \hat{S}_1, \hat{I}_1) - \Phi_1(\hat{Q}_1 + z, \hat{S}_1, \hat{I}_1)\right)}\right\} v(\mathrm{d}z)$$

By (94), we obtain

$$-\psi\left\{\Phi_2\left(\hat{Q}_2, \hat{S}_2, \hat{I}_2\right) - \Phi_2\left(\hat{Q}_2 + z, \hat{S}_2, \hat{I}_2\right)\right\}$$
$$\geq -\psi\left\{\Phi_1\left(\hat{Q}_1, \hat{S}_1, \hat{I}_1\right) - \Phi_1\left(\hat{Q}_1 + z, \hat{S}_1, \hat{I}_1\right)\right\} - \beta\psi\left(\left(\hat{Q}_1 + z\right)^2 - \hat{Q}_1^2\right) - \beta\psi\left(\left(\hat{Q}_2 + z\right)^2 - \hat{Q}_2^2\right), \qquad (105)$$

and thus

$$\frac{\hat{Q}_2 + A}{\psi} \int_\delta^\infty \left\{e^{-\psi\left(\Phi_2(\hat{Q}_2, \hat{S}_2, \hat{I}_2) - \Phi_2(\hat{Q}_2 + z, \hat{S}_2, \hat{I}_2)\right)} - e^{-\psi\left(\Phi_1(\hat{Q}_1, \hat{S}_1, \hat{I}_1) - \Phi_1(\hat{Q}_1 + z, \hat{S}_1, \hat{I}_1)\right)}\right\} v(\mathrm{d}z)$$

$$\geq \frac{\hat{Q}_2 + A}{\psi} \int_\delta^\infty \left\{e^{-\psi\left\{\Phi_1(\hat{Q}_1, \hat{S}_1, \hat{I}_1) - \Phi_1(\hat{Q}_1 + z, \hat{S}_1, \hat{I}_1)\right\} - \beta\psi\left(\left(\hat{Q}_1 + z\right)^2 - \hat{Q}_1^2 + \left(\hat{Q}_2 + z\right)^2 - \hat{Q}_2^2\right)} - e^{-\psi\left(\Phi_1(\hat{Q}_1, \hat{S}_1, \hat{I}_1) - \Phi_1(\hat{Q}_1 + z, \hat{S}_1, \hat{I}_1)\right)}\right\} v(\mathrm{d}z). \qquad (106)$$

$$= -\frac{\hat{Q}_2 + A}{\psi} \int_\delta^\infty e^{-\psi\left(\Phi_1(\hat{Q}_1, \hat{S}_1, \hat{I}_1) - \Phi_1(\hat{Q}_1 + z, \hat{S}_1, \hat{I}_1)\right)} \left\{1 - e^{-\psi\beta\left(\left(\hat{Q}_1 + z\right)^2 - \hat{Q}_1^2 + \left(\hat{Q}_2 + z\right)^2 - \hat{Q}_2^2\right)}\right\} v(\mathrm{d}z)$$

Based on the boundedness assumption of $\Phi_1$, we obtain

$$-\int_\delta^\infty e^{-\psi\left(\Phi_1(\hat{Q}_1, \hat{S}_1, \hat{I}_1) - \Phi_1(\hat{Q}_1 + z, \hat{S}_1, \hat{I}_1)\right)} \left\{1 - e^{-\psi\beta\left(\left(\hat{Q}_1 + z\right)^2 - \hat{Q}_1^2 + \left(\hat{Q}_2 + z\right)^2 - \hat{Q}_2^2\right)}\right\} v(\mathrm{d}z)$$
$$\geq -e^{2\psi C^{(2)}} \int_\delta^\infty \left\{1 - e^{-\psi\beta\left(\left(\hat{Q}_1 + z\right)^2 - \hat{Q}_1^2 + \left(\hat{Q}_2 + z\right)^2 - \hat{Q}_2^2\right)}\right\} v(\mathrm{d}z) \qquad . \qquad (107)$$

with a constant $C^{(2)} > 0$. Because of $1 - e^{-\psi\beta u} \leq \psi\beta u$ for $u \geq 0$, we obtain the estimate

$$-e^{2\psi C^{(2)}} \frac{\hat{Q}_2 + A}{\psi} \int_\delta^\infty \left\{1 - e^{-\psi\beta\left(\left(\hat{Q}_1 + z\right)^2 - \hat{Q}_1^2 + \left(\hat{Q}_2 + z\right)^2 - \hat{Q}_2^2\right)}\right\} v(\mathrm{d}z)$$

$$\geq -e^{2\psi C^{(2)}} \frac{\hat{Q}_2 + A}{\psi} \int_\delta^\infty \psi\beta\left(\left(\hat{Q}_1 + z\right)^2 - \hat{Q}_1^2 + \left(\hat{Q}_2 + z\right)^2 - \hat{Q}_2^2\right) v(\mathrm{d}z) \text{ as } \varepsilon \to +0. \qquad (108)$$

$$\to -2e^{2\psi C^{(2)}} \beta\left(\hat{Q} + A\right) \int_0^\infty \left(\left(\hat{Q} + z\right)^2 - \hat{Q}^2\right) v(\mathrm{d}z)$$

By $\rho - M_1 > 0$, (98), and the boundedness assumption, there is a sufficiently small $\psi > 0$ such that $\rho - e^{2\psi C^{(2)}} M_1 > 0$. For such a small $\psi$, we obtain (99) if the remaining term in (104) is non-negative under the same limit. We show this in what follows. The Hölder continuity and boundedness of $\Phi_1$ give

$$\frac{\hat{Q}_1-\hat{Q}_2}{\psi}\int_\delta^\infty\left\{1-e^{-\psi\left(\Phi_1(\hat{Q}_1,\hat{S}_1,\hat{I}_1)-\Phi_1(\hat{Q}_1+z,\hat{S}_1,\hat{I}_1)\right)}\right\}v(\mathrm{d}z)\geq -\frac{|\hat{Q}_1-\hat{Q}_2|}{\psi}\int_\delta^\infty\left\{e^{\psi C^{(0)}z^\theta}-1\right\}v(\mathrm{d}z). \qquad (109)$$

We assume (if necessary) a smaller $\psi$ with $\psi C^{(0)}<b$. For a small $z>0$, we have $\left\{e^{\psi C^{(0)}z^\theta}-1\right\}v(\mathrm{d}z)\sim z^{\theta-\alpha-1}e^{-bz}\mathrm{d}z$, implying its integrability near $z=0$ by $\theta>\alpha$. For a large $z$, we have $\left\{e^{\psi C^{(0)}z^\theta}-1\right\}v(\mathrm{d}z)\sim az^{\theta-\alpha-1}e^{-(bz-\psi C^{(0)}z^\theta)}\mathrm{d}z$, suggesting far-field integrability. We then conclude that the right-most term of (109) is bounded, and by (69) obtain

$$\frac{\hat{Q}_1-\hat{Q}_2}{\psi}\int_\delta^\infty\left\{1-\exp\left(-\psi\left(\Phi_1(\hat{Q}_1,\hat{S}_1,\hat{I}_1)-\Phi_1(\hat{Q}_1+z,\hat{S}_1,\hat{I}_1)\right)\right)\right\}v(\mathrm{d}z)\to 0 \text{ as } \varepsilon\to +0 \qquad (110)$$

for each $\beta>0$.

Using (104), (108), and (110), we obtain $H_1-H_2\leq 0$ as in the uncertainty-neutral case.

### A.3 Proof of Proposition 3

Firstly, consider the uncertainty-neutral case ($\psi\to +0$). Set $(i_0,j_0,l_0)=\arg\max_{i,j,l}\Phi_{i,j,l}$. We then obtain

$$\Xi_{i_0,j_0,l_0}\Phi_{i_0,j_0,l_0}=\Theta_{i_0,j_0,l_0}[\Phi]-h+f(S_{j_0})\leq\Theta_{i_0,j_0,l_0}\left[\Phi_{i_0,j_0,l_0}\right]-h+f(S_{j_0}), \qquad (111)$$

where $\Theta_{i_0,j_0,l_0}\left[\Phi_{i_0,j_0,l_0}\right]$ is $\Theta_{i_0,j_0,l_0}[\Phi]$ whose arguments are all $\Phi_{i_0,j_0,l_0}$. Hence, it holds that

$$-h+f(S_{j_0})\geq\Xi_{i_0,j_0,l_0}\Phi_{i_0,j_0,l_0}-\Theta_{i_0,j_0,l_0}\left[\Phi_{i_0,j_0,l_0}\right]\geq -\frac{o}{W}, \qquad (112)$$

and we obtain $h\leq f(S_{j_0})+oW^{-1}\leq f(0)+oW^{-1}$. The other bound follows by a similar argument. For the uncertainty-averse case ($\psi>0$), the same bounds hold true because $J_{i,j,l,k}$ in (37) is increasing with respect to $\Phi_{i,j,l}-\Phi_{i+k,j,l}$ and $\Phi_{i,j,l}-\Phi_{i-1,j,l}$, and it vanishes if all $\Phi_{i+k,j,l}$ and $\Phi_{i-1,j,l}$ are replaced by $\Phi_{i,j,l}$.

### A.4 Proof of Proposition 4

First, we introduce several notations used in the proof. Based on the construction, the scheme is degenerate elliptic for $\{\Phi_{i,j,l}\}$ (Definition 2 of Oberman [83]) because (38) is rewritten as

$$h-f(S_j)+\Xi_{i,j,l}\Phi_{i,j,l}-\Theta_{i,j,l}[\Phi]+\Lambda_{i,j,l}[\Phi]\equiv\sum_{(i',j',l')\neq(i,j,l)}\vartheta_{i,j,l,i',j',l'}\left(\Phi_{i,j,l}-\Phi_{i',j',l'}\right)+h-f(S_j)=0. \qquad (113)$$

Here, each $\vartheta_{i,j,l,i',j',l'}:\mathbb{R}\to\mathbb{R}$ is non-decreasing, and for each $i,j,l$ at least one $\vartheta_{i,j,l,i',j',l'}$ is strictly increasing with $\dot{\vartheta}_{i,j,l,i',j',l'}>0$, where the dot "•" means the derivative. We say that the vertices $(i_0,j_0,l_0)$ and $(i_1,j_1,l_1)$ are connected if we have a path such that

$$(i_0, j_0, l_0) = (i_{0,0}, j_{0,0}, l_{0,0}) \to (i_{0,1}, j_{0,1}, l_{0,1}) \to ... \to (i_{0,M}, j_{0,M}, l_{0,M}) = (i_1, j_1, l_1) \tag{114}$$

with some sequence $\{(i_{0,m}, j_{0,m}, l_{0,m})\}_{m=0,1,2,...,M}$ and $M \in \mathbb{N}$ such that $\dot{\vartheta}_{i_{0,m-1}, j_{0,m-1}, l_{0,m-1}, i_{0,m}, j_{0,m}, l_{0,m}}(\cdot) > 0$ ($m = 0,1,2,...,M$).

We assume that the statement of the proposition is false. Then, there is one $(i_0, j_0, l_0)$ with $\Phi_{i_0,j_0,l_0} - \Psi_{i_0,j_0,l_0} = \max_{i,j,l}\{\Phi_{i,j,l} - \Psi_{i,j,l}\} > 0$, such that $\Phi_{i_0,j_0,l_0} - \Phi_{i,j,l} \geq \Psi_{i_0,j_0,l_0} - \Psi_{i,j,l}$ for all $(i,j,l)$ and $\Phi_{i_0,j_0,l_0} - \Phi_{i_1,j_1,l_1} > \Psi_{i_0,j_0,l_0} - \Psi_{i_1,j_1,l_1}$ for some $(i_1,j_1,l_1)(\neq (i_0,j_0,l_0))$. Under this assumption, there is a path connecting $(i_0,j_0,l_0)$ and $(i_1,j_1,l_1)$ having the form (114). Furthermore, there is one $0 \leq \bar{m} \leq M-1$ with

$$\Phi_{i_{0,\bar{m}}, j_{0,\bar{m}}, l_{0,\bar{m}}} - \Phi_{i_{0,\bar{m}+1}, j_{0,\bar{m}+1}, l_{0,\bar{m}+1}} > \Psi_{i_{0,\bar{m}}, j_{0,\bar{m}}, l_{0,\bar{m}}} - \Psi_{i_{0,\bar{m}+1}, j_{0,\bar{m}+1}, l_{0,\bar{m}+1}}. \tag{115}$$

If such an $\bar{m}$ does not exist, then for all $0 \leq m \leq M-1$, we obtain $\Phi_{i_{0,m}, j_{0,m}, l_{0,m}} - \Phi_{i_{0,m+1}, j_{0,m+1}, l_{0,m+1}} \leq \Psi_{i_{0,m}, j_{0,m}, l_{0,m}} - \Psi_{i_{0,m+1}, j_{0,m+1}, l_{0,m+1}}$. Summing up this inequality for all $0 \leq m \leq M-1$ leads to a contradiction $\Phi_{i_0,j_0,l_0} - \Phi_{i_1,j_1,l_1} \leq \Psi_{i_0,j_0,l_0} - \Psi_{i_1,j_1,l_1}$. Using (115), we obtain the strict inequality

$$\begin{aligned}&\vartheta_{i_{0,m}, j_{0,m}, l_{0,m}, i_{0,m+1}, j_{0,m+1}, l_{0,m+1}} \left(\Phi_{i_{0,m}, j_{0,m}, l_{0,m}} - \Phi_{i_{0,m+1}, j_{0,m+1}, l_{0,m+1}}\right) \\ &> \vartheta_{i_{0,m}, j_{0,m}, l_{0,m}, i_{0,m+1}, j_{0,m+1}, l_{0,m+1}} \left(\Psi_{i_{0,m}, j_{0,m}, l_{0,m}} - \Psi_{i_{0,m+1}, j_{0,m+1}, l_{0,m+1}}\right)\end{aligned} \tag{116}$$

and thus

$$f(S_j) - \bar{h} = \sum_{(i',j',l') \neq (i,j,l)} \vartheta_{i,j,l,i',j',l'}(\Phi_{i,j,l} - \Phi_{i',j',l'}) > \sum_{(i',j',l') \neq (i,j,l)} \vartheta_{i,j,l,i',j',l'}(\Psi_{i,j,l} - \Psi_{i',j',l'}) = f(S_j) - \underline{h} \tag{117}$$

at $(i,j,l) = (i_{0,\bar{m}}, j_{0,\bar{m}}, l_{0,\bar{m}})$, which is a contradiction to $\bar{h} \geq \underline{h}$. Therefore, $\{\Phi_{i,j,l}\} \leq \{\Psi_{i,j,l}\}$ in component-wise. The second statement follows immediately by exchanging the role of $\{\Phi_{i,j,l}\}$, $\{\Psi_{i,j,l}\}$.

**Remark 8** For each $l \in \{0,1,2,...,\bar{L}\}$, the existence of paths connecting two arbitrary vertices in the $i$ and $j$ directions follows from (28), (30), and (36). The small assumption of $\bar{L}$ much avoid $\max_{i,j} L^*(Q_i, S_j) < \bar{L}$. In addition, $\bar{S}$ needs to avoid $\max_{i,j}\{\eta^*(Q_i, S_j) + S_j\} < \bar{S}$. **Figure 10** shows a typical connection among the vertices.

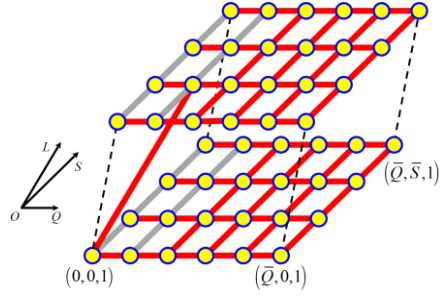

**Figure 10.** A connection among the vertices for $(N_Q, N_S, \bar{L}) = (5,5,2)$. Each red edge represents the connection between its boundary vertices, and the others are colored grey. The connections among the vertices of $l = 1$ and $l = 2$ are all omitted except for one for illustration purposes. The grey edges for a small $Q$ are due to $F = 0$ by $Q \leq \hat{Q}$.

**Appendix B**

The process $(\hat{\phi}_t)_{t \geq 0}$ such that $|\hat{\phi}_t| \leq \bar{\phi} < b$ with a small $\bar{\phi} > 0$ satisfies (7) and (8). For (7), we have

$$
\begin{aligned}
&\mathbb{E}_{\mathbb{P}}\left[\int_0^T \int_0^\infty \int_0^{Q_{s-}+A} (\phi_s(z) - 1 - \ln \phi_s(z)) \, \mathrm{d}u v(\mathrm{d}z) \, \mathrm{d}s\right] \\
&\leq \mathbb{E}_{\mathbb{P}}\left[\int_0^T (Q_{s-} + A) \int_0^\infty (\exp(\bar{\phi} z) - 1 - \bar{\phi} z) v(\mathrm{d}z) \, \mathrm{d}s\right] \\
&= \mathbb{E}_{\mathbb{P}}\left[\int_0^T (Q_{s-} + A) \, \mathrm{d}s\right] \int_0^\infty (\exp(\bar{\phi} z) - 1 - \bar{\phi} z) v(\mathrm{d}z) \\
&< +\infty
\end{aligned}
\tag{118}
$$

For (8), we have

$$
\begin{aligned}
&\mathbb{E}_{\mathbb{P}}\left[\exp\left\{\int_0^T \int_0^\infty \int_0^{Q_{s-}+A} (1 - \phi_s(z) + \phi_s(z) \ln \phi_s(z)) \, \mathrm{d}u v(\mathrm{d}z) \, \mathrm{d}s\right\}\right] \\
&= \mathbb{E}_{\mathbb{P}}\left[\exp\left\{\int_0^\infty (1 - \exp(\bar{\phi} z) + \bar{\phi} z \exp(\bar{\phi} z)) v(\mathrm{d}z) \int_0^T (Q_{s-} + A) \, \mathrm{d}s\right\}\right] \\
&= \mathbb{E}_{\mathbb{P}}\left[\exp\left\{C \int_0^T (Q_{s-} + A) \, \mathrm{d}s\right\}\right] \\
&\leq \exp(CAT) \mathbb{E}_{\mathbb{P}}\left[\int_0^T \exp(CQ_{s-}) \, \mathrm{d}s\right]
\end{aligned}
\tag{119}
$$

with $C \equiv \int_0^\infty (1 - \exp(\bar{\phi} z) + \bar{\phi} z \exp(\bar{\phi} z)) v(\mathrm{d}z) \in (0, \infty)$. This $C$ becomes arbitrarily small when choosing a small $0 < \bar{\phi} < b$. We then apply Theorem 2.16(ii), Lemma 5.3, Lemma 6.5, and Proposition 11.2 of Duffie et al. [100] to $(Q_t)_{t \geq 0}$ and obtain that the last expectation of (119) is bounded because $\int_0^\infty (\exp(\vartheta z) - 1) v(\mathrm{d}z)$ (the non-affine part of the characteristics) is analytic for all complex variables $\vartheta \in U$ with a neighborhood $U$ of the origin in a complex plane, such that all elements of $U$ have real parts smaller than $b$.

**Appendix C**

We derive the autocorrelation function of $(Q_t)_{t \geq 0}$ for a stationary state. The autocorrelation $\omega(\delta)$ for the lag $\delta \geq 0$ is defined as

$$\omega(\delta) = \frac{\mathbb{E}\left[(Q_t - \mathbb{E}[Q_t])(Q_{t+\delta} - \mathbb{E}[Q_t])\right]}{\mathbb{E}\left[(Q_t - \mathbb{E}[Q_t])(Q_t - \mathbb{E}[Q_t])\right]} = \frac{\mathbb{E}\left[Q_t Q_{t+\delta} - \mu^2\right]}{\mathbb{E}\left[Q_t^2 - \mu^2\right]} \tag{120}$$

with the stationary mean $\mu \equiv \mathbb{E}[Q_t] = (\rho - M_1)^{-1}(AM_1 + \rho \underline{Q})$. By (1), Itô's formula gives

$$d(Q_t Q_{t+\delta}) = Q_t \rho (\underline{Q} - Q_{t+\delta}) d\delta + Q_t \Delta Z_{t+\delta}, \tag{121}$$

where $\Delta Z_t$ formally represents the jump generated by $N$ at $t$. The jump at $t+\delta$ is generated by the kernel $Q_{(t+\delta)-} \rho a z^{-(1+\alpha)} e^{-bz} dz$, suggesting that

$$\mathbb{E}[Q_t \Delta Z_{t+\delta}] = \mathbb{E}\left[Q_t (A + Q_{t+\delta}) \int_0^\infty \frac{a}{z^{1+\alpha}} z e^{-bz} dz\right] d\delta = M_1 \left(\mathbb{E}[Q_t Q_{t+\delta}] + A \mathbb{E}[Q_t]\right) d\delta. \tag{122}$$

We then obtain

$$\begin{aligned}\frac{d}{d\delta} \mathbb{E}[Q_t Q_{t+\delta}] &= -\rho \mathbb{E}[Q_t Q_{t+\delta}] + \rho \underline{Q} \mathbb{E}[Q_t] + M_1 \mathbb{E}[Q_t Q_{t+\delta}] + AM_1 \mathbb{E}[Q_t] \\ &= -(\rho - M_1) \mathbb{E}[Q_t Q_{t+\delta}] + (\rho \underline{Q} + AM_1) \mathbb{E}[Q_t]\end{aligned}, \tag{123}$$

$$\frac{d}{d\delta} \mathbb{E}\left[Q_t Q_{t+\delta} - \mu^2\right] = -(\rho - M_1) \mathbb{E}[Q_t Q_{t+\delta}] + (\rho - M_1) \mu \cdot \mu = -(\rho - M_1) \mathbb{E}\left[Q_t Q_{t+\delta} - \mu^2\right], \tag{124}$$

and hence

$$\frac{d}{d\delta} \omega(\delta) = -(\rho - M_1) \omega(\delta) \quad \text{with} \quad \omega(0) = 1 \quad \text{or equivalently} \quad \omega(\delta) = e^{-(\rho - M_1)\delta}. \tag{125}$$

**Appendix D**

We formally derive the moments $m_n(t) = \mathbb{E}[Q_t^n]$ ($n \in \mathbb{N}$): By (1), for $n \in \mathbb{N}$, Itô's formula yields

$$d(Q_t^n) = n Q_t^{n-1} \rho (\underline{Q} - Q_t) dt + (Q_{t-} + \Delta Z_t)^n - Q_{t-}^n = (n\rho \underline{Q} Q_t^{n-1} - n\rho Q_t^n) dt + \sum_{k=0}^{n-1} \binom{n}{k} Q_{t-}^k (\Delta Z_t)^{n-k}. \tag{126}$$

By taking the expectation of (126) with

$$\mathbb{E}\left[Q_{t-}^k (\Delta Z_t)^{n-k}\right] = \mathbb{E}\left[Q_t (A + Q_t) \int_0^\infty \frac{a}{z^{1+\alpha}} z^{n-k} e^{-bz} dz\right] dt = M_{n-k} (m_{k+1} + A m_k) dt, \tag{127}$$

we obtain

$$\frac{dm_n}{dt} = n\rho \underline{Q} m_{n-1} - n\rho m_n + \sum_{k=0}^{n-1} \binom{n}{k} M_{n-k} (m_{k+1} + A m_k). \tag{128}$$

By $\rho - M_1 > 0$, we recursively obtain a stationary $m_n$ ($n = 1, 2, 3, \ldots$) from (128) by omitting the temporal differential term.

**Acknowledgements**

Kurita Water and Environment Foundation 19B018 and 20K004, and Grant for Environmental Research Projects from the Sumitomo Foundation 203160 support this research.

**Declaration**

The authors have no conflict of interests to declare.

**Contributions**

HY: Original manuscript, authorization, conceptualization, methodology, formal analysis, computation, review, and editing
MT: Formal analysis, Review, Editing